\RequirePackage{fix-cm}
\documentclass[11pt,a4paper]{article}   
  
\usepackage{graphicx}
\usepackage[utf8]{inputenc}
\usepackage{amsmath}
\usepackage{mathptmx}
\usepackage{enumerate}
\usepackage[misc]{ifsym}
\usepackage{ArticleEfficient}
 \usepackage[]{algorithm2e}
 \usepackage[compatible]{algpseudocode}
 \usepackage{caption}
 \usepackage{csvsimple}
 \usepackage{subcaption}
 \usepackage{color}
\usepackage{url}
\title{Efficient partitioning and reordering of conforming virtual element discretizations for large scale Discrete Fracture Network flow parallel solvers}
\author{Stefano Berrone $^*$ \and 
	Alice Raeli$^\dagger$ }

\date{}

\begin{document}
\maketitle


\blfootnote{*S. Berrone (\Letter) 
	\em{stefano.berrone@polito.it} 
}
\blfootnote{
	$\dagger$A. Raeli (\Letter) 
	\em{alice.raeli@polito.it}\newline
	Politecnico di Torino \\
	Dipartimento di Scienze Matematiche, \\
	Corso Duca degli Abruzzi 24, 10129 Torino, Italy}

\begin{abstract}
Discrete Fracture Network models are largely used for very large scale geological flow simulations. For this reason numerical methods require an investigation of tools for efficient parallel solutions on High Performance Computing systems. In this paper we discuss and compare several partitioning and reordering strategies, that result to be highly efficient and scalable, overperforming the classical mesh partitioning approach used to partition a conforming mesh among several processes.
\end{abstract}

\section{Introduction}
\label{intro}
The flow in fractured media is a relevant topic in several engineering applications such as aquifers monitoring, disposal and geological storage of nuclear wastes, prevention of accidental dispersion of contaminants, oil and gas enhanced production and many other applications \cite{lei2017use}, \cite{national1996rock}. When the rock matrix of the geological formation has a very small porosity, the contribution of the rock matrix surrounding fractures can have a marginal impact on the flow pattern. Intensity and direction of flow often depend almost uniquely on the distribution of fractures and on their hydraulic properties. \emph{Discrete Fracture Network} models, DFN in the following, simulate transport and flow within fractured material (usually rocks) using discrete computational strategies to approach the real solution \cite{neuman2005trends}, \cite{jaffre2012modeling}. 

In Section~\ref{DFN} we briefly introduce the DFN flow formulation in an impervious rock matrix, following \cite{cacas1990modeling}, \cite{nordqvist1992variable}, \cite{dershowitz1999derivation}, \cite{fidelibus20072d}, \cite{berrone2017non}, modeling fractures as planar polygons. Brief information concerning the mesh creation and the degrees of freedom handling are given in Subsections~\ref{Sec:MMesh} and~\ref{Sec:DofHand} respectively.

Due to the stochastic nature of the DFN model generated starting from probabilistic distribution of fracture position, orientation, size and of hydraulic parameters, several large simulations are needed in order to perform uncertainty quantification for flow quantities of interest \cite{BERRONE2015603}, \cite{berrone2018uncertainty}, \cite{canuto2019uncertainty}, \cite{pieraccini2020uncertainty}. 
As a consequence of the stochastic generation of the networks, DFNs for practical applications are usually very complex. In fact, DFNs count a  large number of fracture intersections with some critical properties, such as very narrow angles or multiple intersection zones.
The generation of a mesh conforming to fracture intersections in a DFN is often a complicated and challenging process, so different numerical approaches possibly circumventing the problem could be applied: standard or mixed Finite Element Methods \cite{vohralik2007mixed}, \cite{hyman2014conforming},  \cite{sentis2017coupling}, hybrid mortar methods \cite{pichot2010mixed}, \cite{benedetto2016hybrid},  optimization methods \cite{mustapha2007new},\cite{berrone2015parallel},  and others \cite{noetinger2012quasi}. This work focuses on a Virtual Element discretization with a conforming polygonal mesh approach \cite{benedetto2016globally} following the strategy presented in \cite{boriodauria} for the generation of the conforming mesh, nevertheless the implementation presented in this work is independent of the element geometry and can be easily extended to other approaches.

The number of mesh cells required for a DFN simulation depends on the number and the size of fractures, the density of the network, the range of scale-lengths generated by the network and on the accuracy of the approximation sought.  Despite the large scale of realistic 3D DFN geological formations, an efficient parallel High Performance Computing approach enables flow and transport simulations. In \cite{berrone2015parallel},\cite{berrone2019parallel} is presented a parallel \emph{master/slaves} approach associated with an PDE (Partial Differential Equation) constrained optimization approach: this choice guarantees scalability and accuracy on the solution requiring a non-standard solver for the linear systems. The method used in this work is based on a more standard conforming discretization approach and proposes a parallelization of the global problem, based on an equitable distribution of the work load, that minimizes communications and may resorts to common solvers and preconditioners for the linear systems.

An efficient parallel DFN simulation necessitate a high-quality partitioning of the degrees of freedom, such that the computations are well-balanced with minimal communications between processes. A multi-level approach is proposed in \cite{ushijima2019multilevel} where the advantages of a DFN-based partitioning are presented in comparison with a classical mesh-partitioning approach.  In Section~\ref{Sec:ParallelPartitioning} different partitioning strategies of the DFN among the processes are presented resorting to different types of DFN based graphs representation; the chosen tool for problem partitioning is the graph partitioning library METIS \cite{metismanual}. To each parallel process is assigned a subset of fractures such that these subsets are disjoint and all the fractures in the DFN are assigned only to one process. Once the DFN is partitioned among processes, in Section~\ref{loadDofs}, we present a method to number the degrees of freedom with respect to the partitioned DFN. The objective of this ordering is to minimize the communications of the iterative method used to solve the linear system by PETSc \cite{balay2019petsc} toolkit.  
The resolution method used in our context is a preconditioned conjugate gradient with a Jacobi preconditioner, other, more efficient, approaches can be used but are not tested in this paper, as we focus on the effect of the MPI parallel communications not on the efficiency of the linear solver.

Although in this work we refer to the specific case of DFNs, the methods presented are compatible with other problems on networks and graphs with a relevant cost of computational operations on nodes and edges; the same approach can be applied to other numerical methods possibly changing the structure of the mesh. In Section~\ref{Sec:NumericalRes} we conclude with numerical tests that investigate the partitioning/reordering methods mentioned above.

\section{The Discrete Fracture Network Discretization}\label{DFN}

In this section, we briefly introduce the notation used in the following. We assume a Darcy flow model inside the fractures and an impervious surrounding rock matrix. Moreover, we assume continuity of the hydraulic head and conservation of the flux at the fracture intersections. We do not provide deeper details of the flow model, as they are not relevant for the focus of the paper and the proposed methods very weakly depend on these modeling choices. Full details concerning the model used for the numerical tests can be found in \cite{berrone2017non} and references therein.

\begin{figure}[!h]
	\begin{center}
		\includegraphics[width=.4\textwidth]{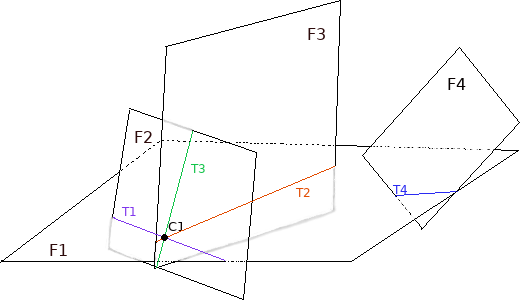}
	\end{center}
	\caption{A simple DFN with four planar rectangular fractures and four intersections.}
	\label{fig:dfn_ex}
\end{figure}

Let $\Omega$ be the DFN, given by the union of planar polygonal fractures $F_r, r=1,...,\#\mathcal{F}$, intersecting each other to form traces $T_m,m=1,...,\#\mathcal{T}$, see Figure~\ref{fig:dfn_ex}.  We refer to fractures and traces as \emph{domains} and \emph{interfaces}. A multiple domains network will represent the DFN $\Omega:=\bigcup_{r=1,...,\#\mathcal{F}}F_r$.

Let $\mathcal{F}$ be the set of the fractures and $\mathcal{T}$ the set of traces; we assume that each trace $T_m$, for $m=1,\ldots, \#\mathcal{T}$, is given by the intersection of two fractures $T_m=F_r\cap F_s$. This assumption induce a map between each trace index and a couple of fracture indexes, $IT(m) = (r,s)$ with $r < s$, such that $F_r\cap F_s=T_m$.  We define $\mathcal{C\!\!\!P}$ the set of cross points of the DFN; we assume that a cross point is a multiple intersection point between three traces, they so belong to three fractures as well. The induced intersection map is defined such that $ICP(t)=(r,s,q)\,$with $r<s<q$, and $F_r\cap F_s\cap F_q=CP_t,\,\ \forall t=1,\ldots\#\mathcal{C\!\!\!P}$.

In Figure~\ref{fig:dfn_ex} an example with four fractures is given, so we have $\mathcal{F}:=\bigcup_{r=1,...,4}F_i$. This network contains four traces and we can list the following induced maps: $IT(1) = (1,2)$, $IT(2) = (1,3)$, $IT(3) = (2,3)$ and $IT(4) = (1,4)$. 
Moreover, the point $CP_1$ given by the intersection of $F_1, F_2$ and $F_3$, such that $ICP(1)=(1,2,3)$, or by the intersection of traces $T_1, T_2$ and $T_3$, is a cross point. 

\begin{figure*}[!h]
	\begin{subfigure}[t]{0.45\textwidth}
		\includegraphics[width=\textwidth]{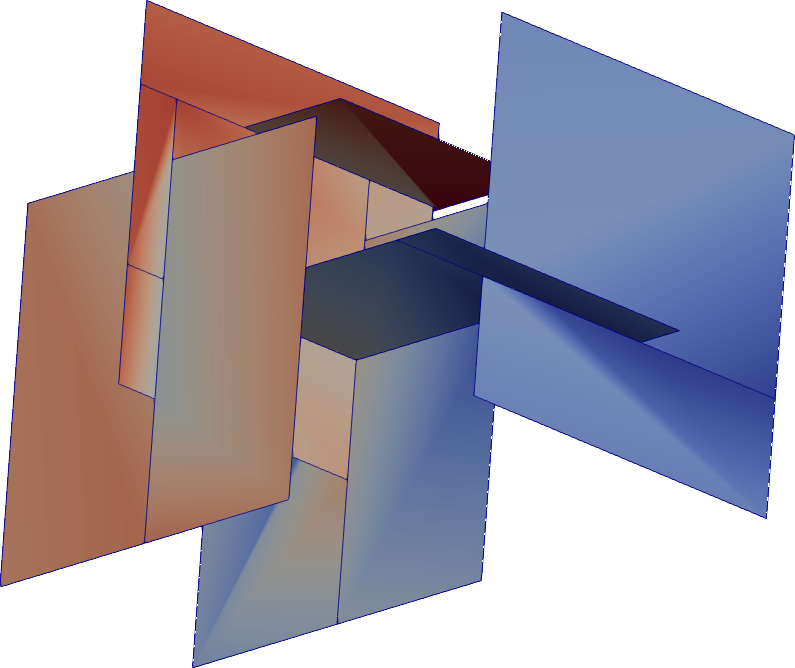}
		\caption{Six Fractures Network }
		\label{fig:F6}
	\end{subfigure}
	~
	\begin{subfigure}[t]{0.45\textwidth}
		\includegraphics[width=\textwidth]{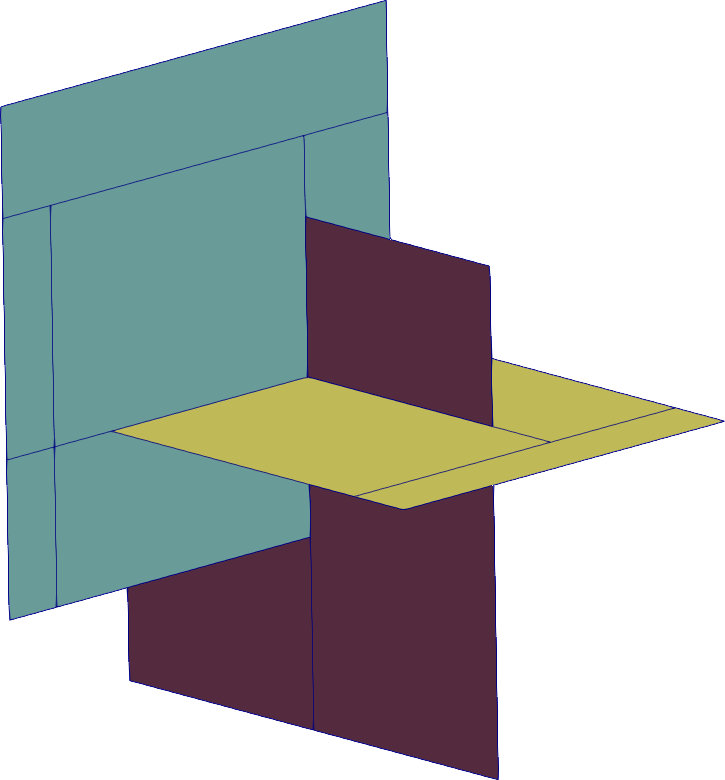}%
		\caption{Cross point, intersection among three fractures. }
		\label{fig:crosspoints}
	\end{subfigure}
	\caption{\emph{Frac6}. A simple DFN with six fractures and six traces.}
	\label{fig:Frac6}
\end{figure*}

In Figure~\ref{fig:Frac6} a DFN with six fractures is represented that will be often used in the following as example. We will refer to it as \emph{Frac6} for the sake of simplicity.

\subsection{The Mesh}\label{Sec:MMesh}
\begin{figure*}[ht!]
	\centering
	\begin{subfigure}[t]{0.45\textwidth}
		\centering
		\includegraphics[height=1.6in]{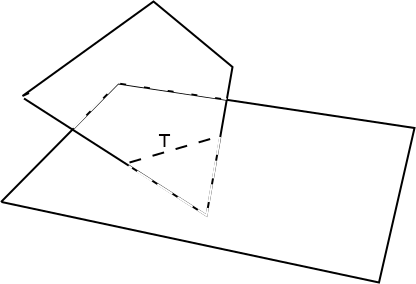}
		\caption{Two intersecting fractures.}
	\end{subfigure}%
	\hfill
	\begin{subfigure}[t]{0.45\textwidth}
		\centering
		\includegraphics[height=1.6in]{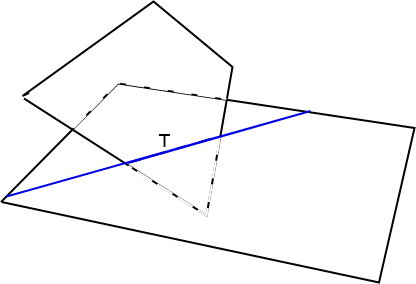}
		\caption{Polygonal mesh on the two fractures induced by trace T.}
		\label{fig:IndByTrace}
	\end{subfigure}
	\caption{Polygonal minimal mesh creation.}
	\label{Fig:InducedMesh}
\end{figure*}

In this section we briefly recall the approach described in \cite{boriodauria} to get a polygonal conforming mesh on a DFN. Let us consider for sake of simplicity two intersecting fractures (Figure~\ref{Fig:InducedMesh}) and let us consider each fracture as a convex cell of a temporary mesh. These cells partially or totally crossed by the trace are split in polygonal convex sub-cells by the trace segments or its possible extensions (see Figure~\ref{fig:F6} and Figure~\ref{fig:IndByTrace}).
This cutting process involves the minimum number of cells on the fracture and is repeated iteratively for all the traces and all the fractures. Full details of this approach can be found in \cite{boriodauria}. The mesh obtained by this process is called \emph{(almost) minimal mesh} (Figure~\ref{Fig:MinimalMesh}).
The minimal mesh can be then refined using, for example, a uniform refinement or an a posteriori mesh refinement as described in \cite{boriodauria}.

In the following we denote by $ E_e\in\mathcal{E},\,e=1,\ldots,\#\mathcal{E}$  a generic cell, and by $N_e$ the set of its nodes.

\subsection{The VEM Discretization}\label{SubSecVEM}
\begin{figure*}[h!]
	\centering
	\begin{subfigure}[t]{.45\textwidth}
		\includegraphics[width=\linewidth]{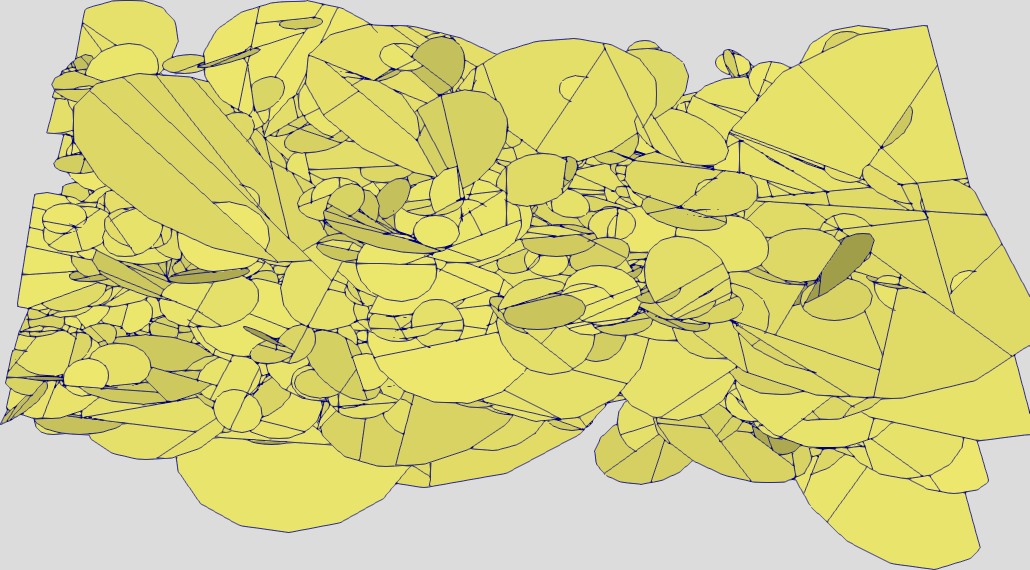}
		\caption{ Minimal mesh induced by traces. }
		\label{Fig:MinimalMesh}
	\end{subfigure}%
	\begin{subfigure}[t]{.45\textwidth}
		\includegraphics[width=\linewidth]{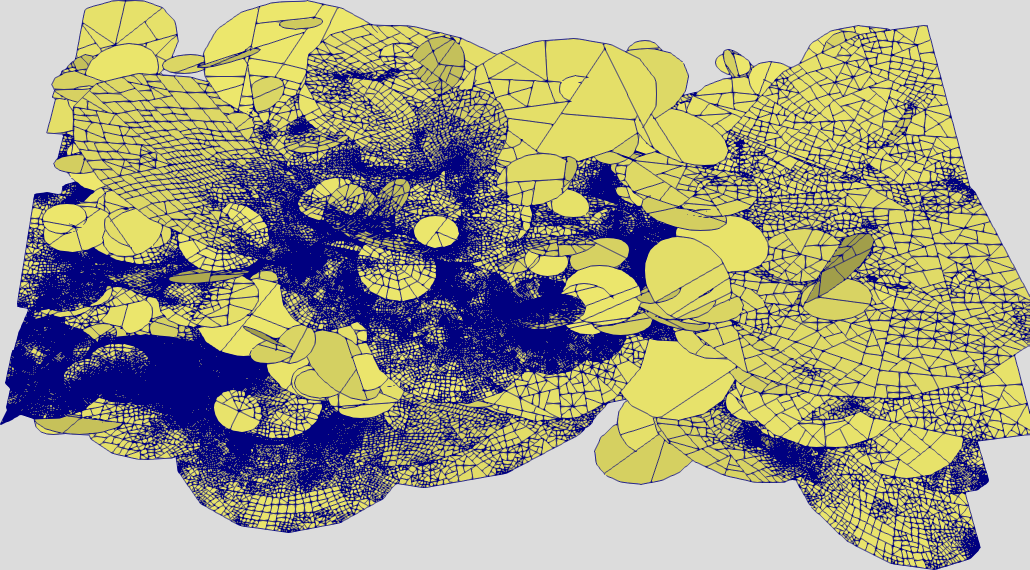}
		\caption{Mesh after 22 adaptive refinements (by momentum cut direction). }
		\label{Fig:22Ref}   
	\end{subfigure}
	\caption{$Frac1000$: Minimal and refined mesh for a VEM resolution. The mesh is refined each iteration by and adaptive mesh refinement algorithm.}
\end{figure*}

The discretization approach applied here is the order $1$ Virtual Element Method that allows the discretization of a second order partial differential equation on a conforming polygonal mesh, \cite{Ahmad2013,Beirao2013a,Beirao2014d,Beirao2015b}. The mesh generation and refinement process here considered yields to convex elements and all the tests proposed are provided on such kind of meshes. Due to the stochastic nature of the DFN the conforming polygonal mesh may contain elements characterized by a low quality due, for example, to a large aspect ratio or the coexistence in the same cell of short and long edges. Although badly shaped elements may affect the quality of a VEM solution, when low polynomial order elements are considered this property is not relevant as shown in \cite{BBorth}.

In order to test the proposed methods on large suitable meshes, in the last Subsection~\ref{Sub:VEM}, we apply the proposed methods on a polygonal mesh obtained applying an {\sl a posteriori} mesh refinement. In this context the problem is solved on progressively refined meshes, being the refinement based on information provided by an {\sl a posteriori} error estimator (Fig.~\ref{Fig:22Ref}).

\subsection{The DOFs Handler}\label{Sec:DofHand}

\begin{figure*}[!ht]
	\begin{subfigure}[t]{0.35\textwidth}
		\centering
		\includegraphics[height=1.5in]{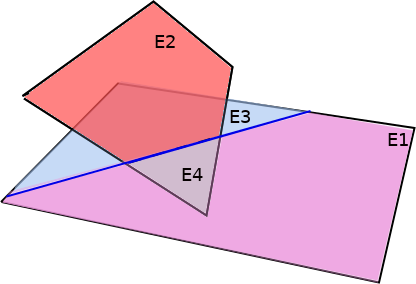}
		\caption{Minimal mesh for two intersecting fractures.}
		\label{fig:IndByTraceCells}
	\end{subfigure}
	\hspace{1.5cm}
	\begin{subfigure}[t]{0.35\textwidth}
		\centering
		\includegraphics[height=1.0in]{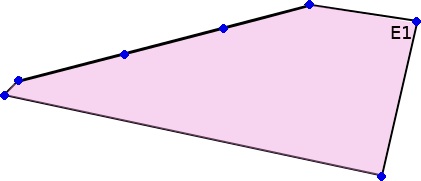}
		\caption{Nodes (blue dots) induced on cell $E_1$.}
		\label{Fig:cellDofs}
	\end{subfigure}
	\caption{Minimal mesh for two intersecting fractures.}
	\label{Figure:buildMesh}
\end{figure*}

A node is a point in space, defined by its coordinates to which we associate a degree of freedom. A degree of freedom ($Dof$) can be defined in many different ways, in the following it will be the point-wise value of the discrete solution on the node considered, and the solution to our problem is written as a linear combination of Lagrangian basis functions $\phi_i$:
\begin{equation}
	\label{uh}
	u_h= \sum_{i=0}^{nDofs-1}Dof_i\phi_i,
\end{equation}
where $\phi_i$ is the basis function related to the $i$-th $Dof$, and $nDofs$ is the total number of degrees of freedom.
In Figure~\ref{fig:IndByTraceCells} is reported the minimal mesh for two intersecting fractures with four cells $E_e$, $e=1,2,3,4$ and in Figure~\ref{Fig:cellDofs} the nodes of the cell $E_1$ are highlighted with blue bullets.

Using a conforming mesh, the degrees of freedom on the DFN include the degrees of freedom at the cross points in $\mathcal{C\!\!\!P}$, the remaining degrees of freedom on the traces in $\mathcal{T}$, and the remaining degrees of freedom inside the fractures in $\mathcal{F}$.
Let $[.]$ be the operator such that $[\mathcal{C\!\!\!P}] $ is the number of degrees of freedom in $\mathcal{C\!\!\!P}$ ($[\mathcal{C\!\!\!P}]=\#\mathcal{C\!\!\!P}$), $[\mathcal{T}]$ is the number of degrees of freedom of the traces (including the degrees of freedom at the cross points), and $[\mathcal{F}]=nDofs$ is the number of degrees of freedom on the fractures (total number of degrees of freedom).

In the following each degree of freedom is uniquely related to an integer index $i=0, \ldots, nDofs-1$ and we refer to a possible permutation of this index set as a reordering. For sake of simplicity, let us assume we are dealing with a problem with homogeneous Dirichlet boundary conditions.
For each node $\mathbf{x}$, not on the boundary of the DFN, let $i$ be the corresponding $Dof$ index, $\phi_i$ the corresponding Lagrangian basis function and $\mathcal{E}_\mathbf{x}=\{E\in \mathcal{E}\, |\, E\cap \mathbf{x} \neq \emptyset \}$ its support, i.e. the set of cells intersecting the node $\mathbf{x}$; we define $Neigh_i$ the set of the Dofs indices corresponding to the vertices of the cells in $\mathcal{E}_\mathbf{x}$ called \emph{neighborhood} of the $i$-th $Dof$. 
Most of the Dofs on fractures have a neighborhood within the fracture, but Dofs on traces and cross points have adjacent cells on different fractures. The neighborhoods of the degrees of freedom of traces and cross points are relevant during the partitioning strategies presented in the following, and, when the connected fractures lie on different MPI processes, they are related to the part of the solution that requires communications between processes.

\section{Parallel Partitioning}\label{Sec:ParallelPartitioning}

\begin{figure*}
	\begin{subfigure}[t]{0.45\textwidth}
		\centering{%
			\includegraphics[width=\textwidth]{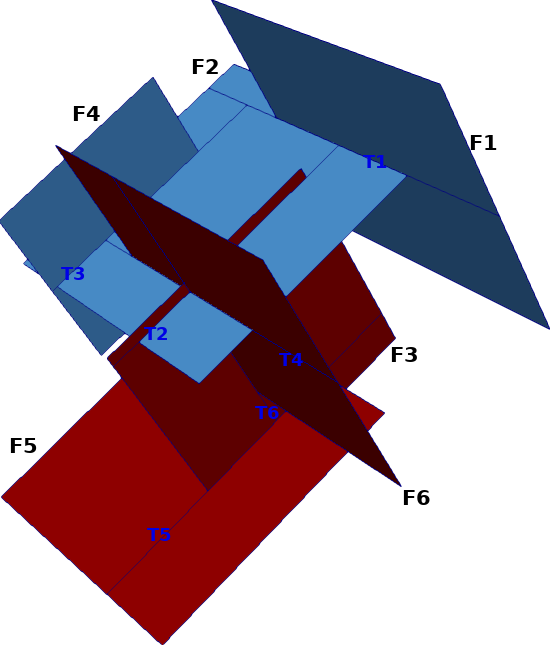}%
		}
		\caption{\emph{Frac6} partition between 2 processes. \\Red fractures are assigned to $P_1$,\\ blue ones are assigned to $P_2$. }
		\label{fig:Frac6Partit}
	\end{subfigure}%
	~ \hspace{1cm}
	\begin{subfigure}[t]{0.4\textwidth}
		\centering{%
			\includegraphics[width=1.1\textwidth]{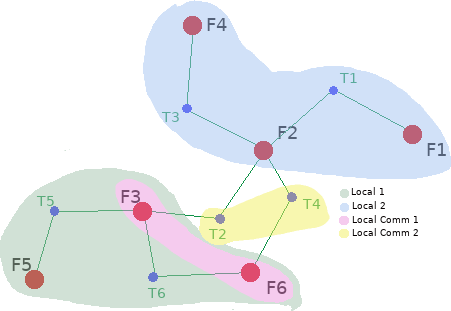}%
		}
		\caption{Graph partitioning among 2 processes. Local and local communicating networks}
		\label{fig:partit}
	\end{subfigure}%
	\caption{\emph{Frac6} - 2 processes partition.}
	\label{Figure:Frac62Procs}
\end{figure*}

The discrete structure of a DFN (union of polygonal fractures) naturally implies several graph representations, \cite{hendrickson1995multi},  \cite{metismanual}, \cite{ushijima2019multilevel}.
In order to balance the computations in DFN flow simulations two main partitioning strategies are possible: the mesh-based partitioning, i.e. the partitioning of the Dofs of the full DFN based on the mesh connectivity, and the DFN-based ones, i.e. the partitioning of the geometrical objects (fractures, traces, cross points) and the corresponding Dofs based on objects connectivity.

In the following we present three different DFN-based graph partitioning strategies among all MPI processes ($\rm{Prs}$). 
Let $\mathcal{P}=\cup_{i=1}^{\#\rm{Prs}} P_i$ be  a partition of fractures, where $P_i$ is the set of fracture assigned to $i$-th process of our parallel environment with $P_i\cap P_j=\emptyset, \forall i\neq j$.
Once the DFN fractures are divided among the processes, auxiliary sub-networks are created in order to suitably manage the communications between processes.

The \textbf{local network} $\mathcal{L}_i$ of the $i$-th process is an ordered set containing the set of fractures $P_i$ and the set of internal traces  $\mathcal{T}_{P_{i}}=\{T_m\in \mathcal{T} \, |\, IT(m)=(r,s),\, F_r,\, F_s \in P_i \}$, i.e. $\mathcal{L}_i=(P_i,\mathcal{T}_{P_i})$.

The \textbf{local communicating network} $\mathcal{LC}_i$ is an ordered set containing a set of fractures and a set of traces that are involved in communications during resolution. Let $\mathcal{F}_{LC,P_i}=\{F_s\in P_i \, |\,\exists T_m\in \mathcal{T} : \, IT(m)=(r,s),\, F_r \notin P_i \}$  be the fractures in $P_i$ that share a trace with a fracture not belonging to $P_i$ having a fracture index smaller than the one of the fracture in $P_i$. Moreover, let $\mathcal{T}_{LC,P_i}=\{T_m\in \mathcal{T} \, |\, IT(m)=(r,s),\, F_r \in P_i, F_s\notin P_i \}$ be the set of traces shared by a fracture of $P_i$ and a fracture not belonging to $P_i$ with a fracture index larger than the one in $P_i$.  We define $\mathcal{LC}_i=(\mathcal{F}_{LC,P_i},\mathcal{T}_{LC,P_i})$. Furthermore, the process that owns a trace in its local communicating network is the one that handles its corresponding Dofs.

We explore in detail the \emph{Frac6} example in Figure~\ref{fig:Frac6Partit} for a two processes partition: red highlighted fractures lie on $\mathcal{L}_1$, light blue ones belong to $\mathcal{L}_2$. In Figure~\ref{fig:partit}  $\mathcal{L}_1$ and  $\mathcal{L}_2$ are disjointed in their graph representation (highlighted two disjointed zones). The traces $T_5$ and $T_6$ lie on process $1$ and they are totally local, so they belong to $\mathcal{L}_1$. Analogously $T_1,\,T_3\in \mathcal{L}_2$.  
The traces $T_2$ and $T_4$ are not local as long they rely two fractures in different processes. Following the construction strategy ($IT(2)=(2,3),IT(4)=(2,6)$) in both cases the two traces are handled by process $2$: their Dofs are indexed with the Dofs of fracture $F_2$, whereas the process $1$ set up the two fractures $F_3$ and $F_6$ to receive Dof indices on these traces from the process 2 (pink subset in Figure~\ref{fig:partit}). We have $\mathcal{LC}_1=\{F_3,\,F_6\},\,$ and $\mathcal{LC}_2=\{T_2,\,T_4\}$.

We define the \textbf{cut} $C$ of a network partition $\mathcal{P}$ as the number of traces that connect fractures belonging to different processes: 
\begin{equation}\label{eq:defCut}
	C=  |\mathcal{C}|,  \quad \mathcal{C}=\{T_m\in \mathcal{T} | IT(m)=(r,s)
	,\, F_r \in P_i, F_s \in P_j, i\neq j \}.
\end{equation}
We note that $C$ is the sum of the number of all the traces contained in local communicating networks of the partitioning.

Let $D_i$ be the number of degrees of freedom of the sub-network associated to the partition $P_i\in\mathcal{P}$, we define the \textbf{imbalance} of a partition $\mathcal{P}$ as:
\begin{equation}\label{eq:defImbal}
	I= \frac{\min_{P_i\in\mathcal{P}} D_i}{\max_{P_i\in\mathcal{P}} D_i} = \frac{D_{min}}{D_{max}}.
\end{equation}
Closer the quantity $I$ is to one, more balanced is the partitioning.

The METIS library \cite{metismanual} used to partition the graph requires the graph stored in its {\em adjacency format}.  
METIS toolkit also requires a {\em partitioning objective} to be chosen between {\em edge-cut minimization} and {\em communication volume minimization}.
By adding {\em weights} to the nodes of a graph we quantify the amount of the computations entrusted to each process of the partition managing those nodes also if an edge-cut minimization is applied, whereas the weights on edges quantify the amount of communication between processes.

\subsection{Partitioning Strategies}\label{Section:Partitioning}

In this section we focus on the partition of the DFN among the computational processes aiming at maximizing data locality in order to minimize communications and at balancing the computational load among the processes. 
Six partitioning strategies are presented and analyzed in the following.

\subsubsection{Partition Graph $Pg$  and Partition Weighted Graph $Wg$}\label{Sub:Pg}

\begin{figure}[!h]
	\centering{%
		\includegraphics[width=.5\textwidth]{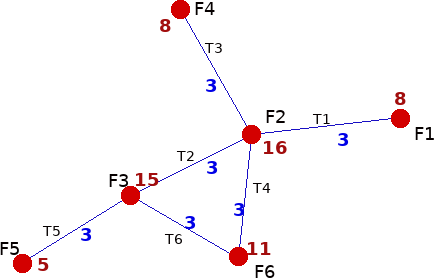}
	}
	\caption{\emph{Frac6}, adjacency graph representation. F stands for fractures and T for traces.
		Red numbers are the weights on nodes and blue ones the weights on edges.}
	\label{fig:Frac6Graph}
\end{figure}

In the simplest graph representation of a DFN the fractures are graph-nodes and the traces are graph-edges (see Figure~\ref{fig:Frac6Graph}), \cite{berrone2015parallel}, \cite{ushijima2019multilevel}; 
weights on nodes and edges are set to $1$. We denote this partitioning as $Pg$ (Partition Graph); this approach aims at minimizing the number of the cut traces (cut-minimization), so the number of communicating traces, but the chosen traces to be cut does not necessarily minimize the amount of data communication during the resolution process. Moreover, we do not consider in the partitioning the amount of computations required by each process.

In order to avoid this problem we can resort to a \emph{weighted} partition ($Wg$) of the graph attaching weights to nodes and edges in the following way:
\begin{itemize}
	\item $\forall F_r\in \mathcal{F}$ the associated node weight is $[F_r]$;
	\item $\forall T_m\in \mathcal{T}$ the associated edge weight is $[T_m]$.
\end{itemize}
By an edge-cut minimization this weighted graph we limit the amount of data communication, whereas the node-weights balance the workload of the processes.

We will present in numerical results that the partition time of the weighted version is almost preserved with respect to the non-weighted version, but the imbalance value $I$ is significantly improved in front of a negligible deterioration of the cut $C$.
\begin{figure*}[!h]
	\centering{%
		\begin{subfigure}[t]{0.45\textwidth}
			\includegraphics[width=\textwidth]{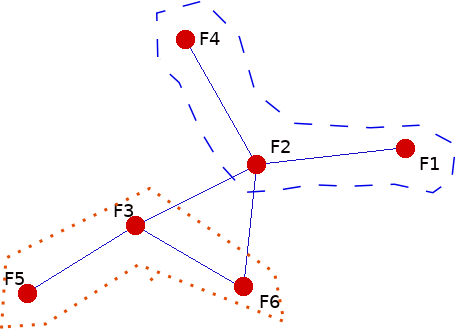}
			\caption{$Pg$ cut minimization partitioning strategy among two processes.}
		\end{subfigure}
		\begin{subfigure}[t]{0.45\textwidth}
			\includegraphics[width=\textwidth]{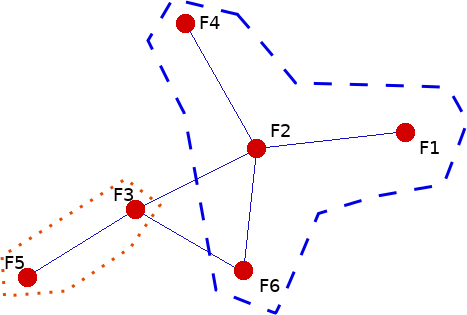}
			\caption{$Wg$ weighted partitioning strategy among two processes.}
		\end{subfigure}
	}
	\caption{\emph{Frac6} partitioning strategies comparison. Nodes in the dotted region lie in $P_1$, dashed ones in $P_2$.}
	\label{fig:Frac6GraphG}
\end{figure*}

In Figure~\ref{fig:Frac6GraphG} an example of the two partitions for $Frac6$ is presented. We note differences on local networks $\mathcal{L}_{0}$ and $\mathcal{L}_{1}$. Increasing the number of the fractures the weighting effect becomes more evident.

\subsubsection{Partition Bipartite Graph $Pb$ and Partition Bipartite Weighted Graph $Wb$}\label{Sub:Pb}
\begin{figure}[!h]
	\centering{%
		\includegraphics[width=.5\textwidth]{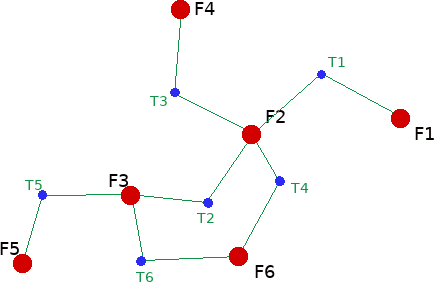}%
	}
	\caption{\emph{Frac6} bipartite graph representation. F stands for fractures and T for traces.}
	\label{fig:Frac6TraceNode}
\end{figure}

Most of the amount of data communication between processes depend on the number of degrees of freedom on the traces in $\mathcal{C}$. In order to increase the impact of the degrees of freedom on the traces on the partitioning process, we introduce a bipartite graph representing the DFN.
In the bipartite graph there are two sets of nodes, the first set representing the fractures and the second one representing the traces. Each element of the first set is only connected with elements of the second one and other way round. The edges correspond to connections between traces and fractures.

When a cut occurs on an edge connecting a trace $T_m$, with $IT(m)=(r,s)$, and a fracture $F_s$, the trace $T_m$ is associated to $\mathcal{LC}_i$ and the fracture $F_s$ is associated to $\mathcal{LC}_j\,$with $i\neq j$. 
In Figure~\ref{fig:Frac6TraceNode} we provide the bipartite graph corresponding to the DFN \emph{Frac6}. 
We denote this partitioning as $Pb$ (Partition Bipartite).
As before we also explore the advantages of a weighted version when we provide weights on nodes and edges ($Wb$ Weighted Bipartite). We assume that the number of degrees of freedom $[T_m]$ on traces is proportional to the number of connections between the trace and the fractures which intersects. We set the edges weights as $[T_m],\, \forall T_m\in\mathcal{T}$. 
In this approach we increase the possibility of the partitioning to entrust traces to different processes and in the weighted version we can distinguish the amount of communication related to the connectivity of the Dofs on the trace and the two connected fractures that can be very different when the size of the elements on the two fractures is different. This may happen, for example, when the two fractures have a strong gap in the transmissivity.
\begin{figure*}[!h]
	\centering{%
		\includegraphics[width=.49\textwidth]{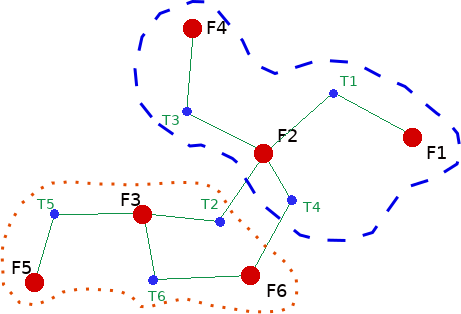}
		\includegraphics[width=.49\textwidth]{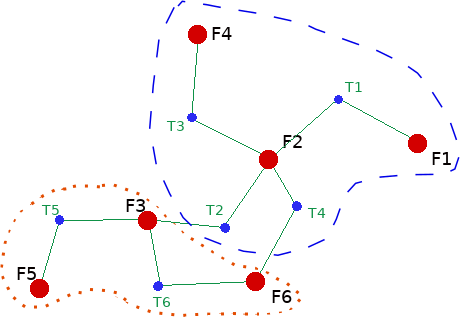}
	}
	\caption{\emph{Frac6}, $Pb$ (left) and $Wb$ (right) partitionings. Nodes in the dotted region are assigned to $P_1$, dashed ones are assigned to $P_2$.}
	\label{fig:Frac6GraphB}
\end{figure*}

In Figure~\ref{fig:Frac6GraphB} $Pb$ and $Wb$  partitionings are represented for the $Frac6$ example. We observe that both strategies for this very small test produce the same local networks $\mathcal{L}_{0}$ and $\mathcal{L}_1$, however the METIS assignment of the nodes differs on $T_2$.

\subsubsection{Partition Tripartite Graph $Pt$  and Partition Tripartite Weighted Graph $Wt$}\label{Sub:pt}

\begin{figure*}[!h]
	\centering{%
		\includegraphics[width=.5\textwidth]{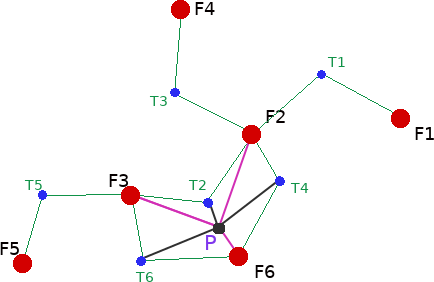}%
	}
	\caption{\emph{Frac6}. Tripartite graph. F denotes fractures, T stands for traces and P corresponds to cross points with multiple intersection. The adjacency matrix considers as connections the green, pink and black edges.}
	\label{fig:Frac6Complete}
\end{figure*}

In this approach we set as nodes of the graph the fractures, the traces and the cross points: they are respectively managed as two, one and zero dimensional domains.  We call tripartite graph the new induced graph on which the edges correspond to connections between traces and cross points, traces and fractures, cross points and fractures.
This approach has a larger number of nodes and further increase the possibility to distribute nodes among processes, moreover the number of edges generated by this graph representation is larger and can change the Dofs distribution among processes if a partition that minimizes the edges-cut is applied.
We notice on the example given in Figure~\ref{fig:Frac6Complete} that a cross point is connected with fractures and traces. 
We denote this partitioning as $Pt$ (Partition Tripartite) in the following. 

Its weighted counterpart $Wt$ (Weighted Tripartite) is built such that:
\begin{itemize}
	\item $\forall F_r\in \mathcal{F}$ the associated node weight is $[F_r]$;
	\item $\forall T_m\in \mathcal{T}$ the corresponding node weight is $[T_m]$; 
	\item $\forall CP_t\in \mathcal{C\!\!\!P}$ the associated node  weights are set to $1$; 
	\item for each edge of the graph connecting a trace to a fracture the corresponding weight is $[T_m]$, the number of Dofs on the trace;
	\item for each edge of the graph connecting a cross point to a fracture, or to a trace, the associated edge weight is $[CP_t]*degree(CP_t)$; $[CP_t]=1$ is the number of Dofs on cross point, and $degree(CP_t)$ is the degree of the node $CP_t$ in the graph; for example in Figure~\ref{fig:Frac6Complete} $degree(CP_t)=6$ as it intersects three fractures and three traces of the DFN.
\end{itemize}

\begin{figure*}[!h]
	\centering{%
		\includegraphics[width=.49\textwidth]{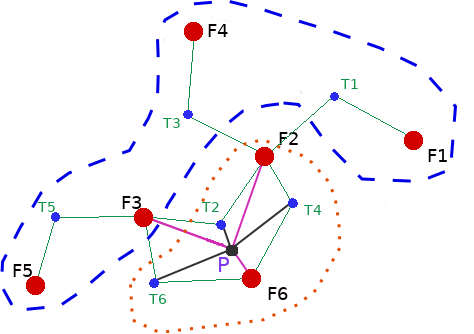}
		\includegraphics[width=.49\textwidth]{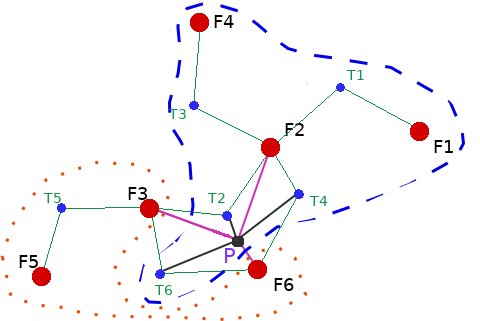}
	}
	\caption{\emph{Frac6}, $Pt$ (left) and $Wt$ (right) partitionings. Nodes in the dotted region lie in $P_1$, dashed ones in $P_2$.}
	\label{fig:Frac6GraphT}
\end{figure*}

In Figure~\ref{fig:Frac6GraphT} we present the $Frac6$ partitioning example, the unweighted and weighted versions differ in their local networks $\mathcal{L}_0$ and $\mathcal{L}_{1}$.

\section{Indexing Dofs}\label{loadDofs}

In this section we focus on global indexing strategies of the degrees of freedom, unique among all the processes. We present a first serial strategy to enumerate the Dofs of the DFN in Algorithm~\ref{Alg:serialDof}. This strategy is among the simplest to be applied to a DFN, and by changing slightly this approach we provide a parallel efficient global indexing.  These algorithms first set indices on domain interfaces of increasing geometrical dimension, then they fill the remaining degrees of freedom with a simple incremental order on fractures.

\begin{algorithm}[]
		\caption{Degrees of Freedom: Serial Assignment.}
	\begin{algorithmic}[1]
		\STATE \ForEach{$ CP_t\in \mathcal{C\!\!\!P}$}
		{
			Assign dof index $d$. $d++$	
		} 
		\STATE \ForEach{$T_m\in \mathcal{T}$}
		{		\ForEach{ node on $T_m$ without index}{
				
				Assign dof index to $d$. $d++$			
			}
		}		
		\STATE \ForEach{$F_r\in \mathcal{F}$}
		{
			\ForEach{node on $F_r$ without index}{
				Assign dof index to $d$. $d++$
			}
		}		
	\end{algorithmic}
		\label{Alg:serialDof}
\end{algorithm} 

\begin{figure*}[!h]
	\begin{subfigure}{0.5\textwidth}
		\includegraphics[width=.95\textwidth]{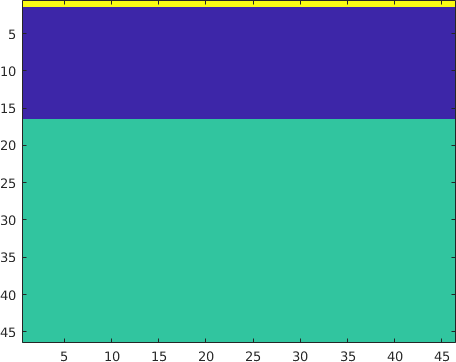}
		\caption{Outline structure: the yellow line refers to cross point test functions, the blue stripe to traces and the green one to fractures internal Dofs.}
	\end{subfigure}
	\begin{subfigure}{0.5\textwidth}
	\includegraphics[width=.8\textwidth]{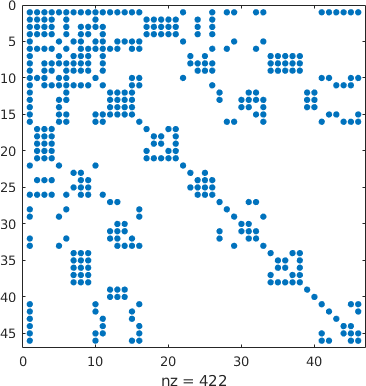}
	\caption{Matrix sparsity pattern.}
	\label{Frac6Spy}
	\end{subfigure}

	\caption{$Frac6, 46$ degrees of freedom distribution for regular DofHandler. }
	\label{Fig:RegularDH}
\end{figure*}

Given $[\mathcal{T}]-[\mathcal{C\!\!\!P}]= M$ the number of degrees of freedom on traces that are not connected to cross points and $[\mathcal{F}]-M-[\mathcal{C\!\!\!P}]=K$ the internal degrees of freedom on fractures; the resulting sparse matrix has the first $[\mathcal{C\!\!\!P}]$ rows concerning the cross points interactions, then $M$ rows concerning the remaining traces Dofs and the last $K$ rows concerning the remaining fracture Dofs. This strategy agglomerates the most communicating rows in the higher part of the matrix, see Figure~\ref{Fig:RegularDH}. This approach is not convenient for a PETSc parallelization (unless different restrictions are imposed) because PETSc subdivides the matrix among the processes in horizontal contiguous stripes.
In order to perform the matrix-vector products needed by the PCG method each process requires updating components of the vector that were computed by other processes at the previous iteration. The elements of the vector that are involved in communications are the elements whose indices are outside the diagonal block of the stripe assigned to the process. In Figure~\ref{Frac6Spy} the sparsity pattern of the matrix highlights the presence of many nonzero elements outside the diagonal blocks for all the processes.
The matrix partition resulting in this case overloads the first process of data communications with almost all the other processes highlighting the negative effect of the latency.
The SpeedUp of Krylov subspace methods depends on the global synchronization during the matrix-vector product communication  \cite{saad2003iterative}, \cite{GMRESReorder}, \cite{kaya2013analysis}.

\subsection{The Reordered DofHandler}
\label{reorderededdofhandler}

We propose a reordering for degrees of freedom in order to enumerate consecutively the Dofs of the geometrical objects contained in the local network of each process. Moreover the degrees of freedom on local communicating traces (including Dofs at cross points) are handled by the process that owns them.
\begin{figure*}[!h]
	\begin{subfigure}{0.5\textwidth}
		\includegraphics[width=.95\textwidth]{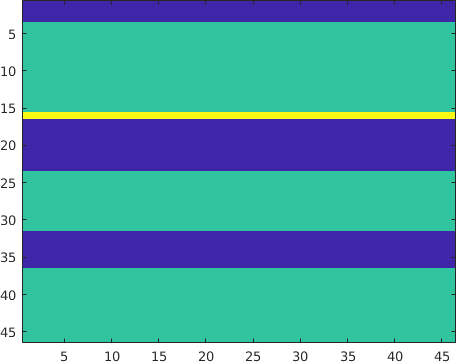}
			\caption{Outline structure: the yellow line refers to cross point test functions, the blue stripe to traces and the green one to fractures  internal Dofs.}
			\label{F6OutlineReo}
	\end{subfigure}
\begin{subfigure}{0.5\textwidth}
\includegraphics[width=.8\textwidth]{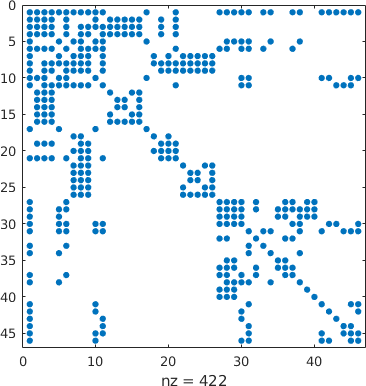}
\caption{Matrix sparsity pattern.}
\label{SpyF6Reord}
\end{subfigure}
	\caption{$Frac6, 46$ degrees of freedom distribution for reordered DofHandler.  3 processes partition.}
	\label{Fig:ReordMat}
\end{figure*}

In Algorithm~\ref{Alg:parDof}, each process numbers its Dofs with a local numbering, then the local numberings are concatenated in the global indexing.
Interfaces indices at traces and cross points are set by the process owning the lowest index fracture and are skipped by the other processes that will inherits the global numbering later. Recalling that $\forall CP_t\in F_r\cap F_s \cap F_q$ the intersection map $ICP(t)=(r,s,q)$ is such that $r<s<q$, and for a trace $T_m\in\, F_r\cap F_s$ we have $IT(m)=(r,s)$ with $r<s$, the algorithm can be sketched in the following way:

\begin{algorithm}[H]
	
	\KwData{Each Process $i$ Call}
	\ForEach{$CP_t\in \mathcal{C\!\!\!P}$, compute $ICP(t)$}{
		\uIf{ $F_r\in \mathcal{L}_i$}{
			Assign dof index to $d_{local}$; $d_{local}++$;
		}
	}
	\ForEach{$T_m\in \mathcal{T}$, compute $IT(m)$}{
		\uIf{ $F_r\in \mathcal{L}_i$}{
			\ForEach{ node on $T_m$ without index}{
				
				Assign dof index to $d_{local}$; $d_{local}++$;
			}
		}
	}
	
	\ForEach{$F_k\in\mathcal{L}_i$}{
		\ForEach{ node on $F_k$ without index}	
		{
			\If{ $F_k\in \mathcal{LC}_i$}{
				Prepare $F_k$ to receive.
			}
			Assign dof index to $d_{local}$; $d_{local}++$;
		}
	}
	\caption{Pre-communicating phase, degrees of freedom local assignment.}
	\label{Alg:parDof}
\end{algorithm}

The global indices are then computed from local ones adding to them an offset corresponding to the sum of indices counted on the previous processes. Each process determines the global indices for the degrees of freedom previously set, then indices of Dofs not managed by the process will be received by a different one containing the interface in the local communicating network.

This approach ensures a communicating part of the matrix to each process, balancing the communications required by the solver at each iteration. 

In Figure~\ref{Fig:ReordMat} we provide a representation of the matrix previously seen in Figure~\ref{Fig:RegularDH} partitioned among 3 processes by this algorithm, the communicating data (Figure~\ref{F6OutlineReo}) are sketched in the first part of each horizontal stripe. Moreover, in Figure~\ref{SpyF6Reord}, the sparsity pattern is presented; we remark that the $Frac6$ is a small test on which the advantages of partitioning and reordering strategies are not evident, however the third processes highlights a larger diagonal local block. Further details are presented for a larger DFN in Section~\ref{Sec:NumericalRes}.

\section{Numerical Results}\label{Sec:NumericalRes}

\begin{figure*}[!ht]
	\begin{subfigure}[t]{.5\textwidth}
		\centering
		\includegraphics[width=\textwidth]{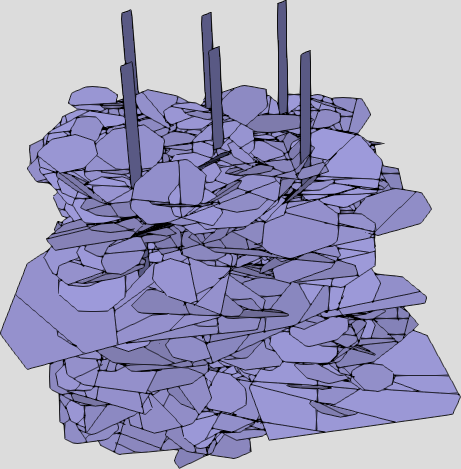}
		\caption{ $Frac512$: Minimal mesh. }
		
	\end{subfigure}%
	\begin{subfigure}[t]{.5\textwidth}
		\centering
		\includegraphics[width=\textwidth]{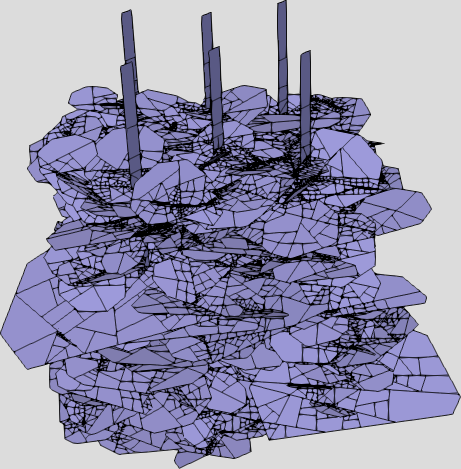}
		\caption{$Frac512$: refined mesh $n$=100. }
		
	\end{subfigure}
	\caption{$Frac512$: initial minimal mesh (left) uniformly refined  until the number of degrees of freedom reaches $[\mathcal{F}]\geq n*512$ (right). }
	\label{Fig:Frac395meshes}
\end{figure*}
\begin{table}[!h]
		\caption{DFN Numerical notations.}
	\label{Notations2}
	\begin{center}
		\small{
			\begin{tabular}{l|l}%
				\hline
				$[\mathcal{F}]$  & Total number of degrees of freedom \\
				\hline
				$\#\mathcal{S}$  & Cardinality of the set $\mathcal{S}$	\\				
				\hline
				$n$  &  Average number of degrees of freedom required on each fracture\\
				\hline
				$C$  & Cut\\
				\hline
				$I$  & Imbalance \\
				\hline
				$Part. Time(s)$  & Partitioning time, in seconds \\
				\hline
				$Res. Time(s)$  & Solution time, in seconds
				\\\hline
			\end{tabular}
		}
	\end{center}
\end{table}
\begin{table}[!h]
		\caption{DFN Partitioning strategies notations.}
	\label{Notations}
	\begin{center}
		\small{
			\begin{tabular}{l|l}%
				\hline
				$Pg$  & Partitioning induced graph. The weights on nodes and edges are $1$. \\ & Fractures as nodes, traces as edges\\
				\hline
				$Wg$  & Weighted partitioning induced graph \\ & Fractures as nodes, traces as edges. \\ & The weights are $[F_r]$ and $[T_m]$ on the associated nodes and edges. \\
				\hline
				$Pb$  & Partitioning bipartite graph. The weights on nodes and edges are $1$. \\ & Fractures and traces as nodes, connections between fractures and traces are  edges\\
				\hline
				$Wb$  & Weighted partitioning bipartite graph \\ & Fractures and traces as nodes, connections between fractures and traces are edges. \\ & The weights are $[F_r]$ and $[T_m]$ on the nodes associated to fractures and traces respectively,\\& $[T_m]$ on the edges connecting traces to fractures. \\
				\hline
				$Pt$  & Partitioning tripartite graph. The weights on nodes and edges are $1$. \\ & Fractures, traces and cross points as nodes, connections between these objects are edges\\
				\hline
				$Wt$  & Weighted partitioning tripartite graph \\ & Fractures, traces and cross points as nodes, connections between these objects are edges. \\ & The weights are:  $[F_r]$ on fracture nodes;\\& 
				$[T_m]$ on trace nodes and on edges connecting trace nodes and fracture nodes; \\ & $6$ on edges connecting cross points to fractures and traces; \\
				&$1$ on nodes representing cross points. 
				\\\hline
			\end{tabular}
		}
	\end{center}
\end{table}

The tests presented in this section concern DFNs with different number of fractures: 512, 1000, 2000 and 4000. We denote the corresponding DFNs as  $Frac512$, $Frac1000$, and so on.
The linear systems of the tests presented are solved using the PETSc preconditioned conjugate gradient iterative method (\emph{KSPCG}) with a Jacobi diagonal preconditioner.  This choice exploits the symmetry and coercivity of the Darcy problem and the corresponding symmetric positive definiteness of the VEM discretization matrix, nevertheless we remark that the partitioning and renumbering methods discussed do not rely on this properties.

\subsection{Reordering Analysis}\label{Sec:ReoAnalys}

\begin{figure*}[h]
	\begin{subfigure}[t]{.45\textwidth}
		\centering
		\includegraphics[width=\textwidth]{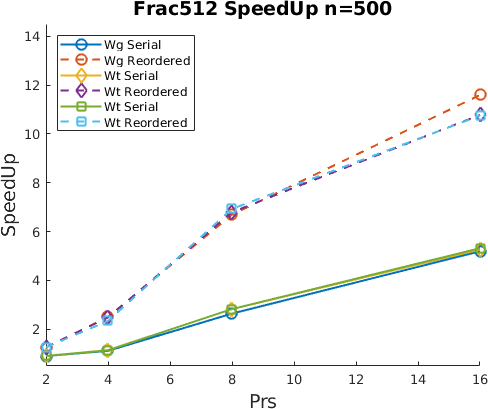}
		\caption{$Frac512\, n=500$: SpeedUp. }
		
	\end{subfigure}
	\begin{subfigure}[t]{.45\textwidth}
		\centering
		\includegraphics[width=\textwidth]{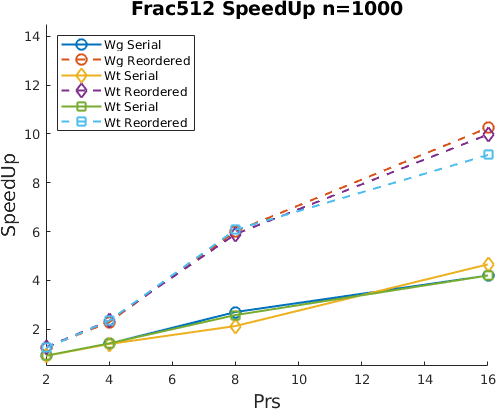}
		\caption{$Frac512\, n=1000$: SpeedUp. }
		
	\end{subfigure}
	\caption{SpeedUp: $Frac512$, weighted partitioning strategies base, bipartite and tripartite. Comparison of reordered and serial DofHandler numbering strategies. $n=500$, left, $n=1000$, right.}
	\label{Fig:Frac512times}
\end{figure*}

\begin{figure*}
	\begin{subfigure}[t]{.45\textwidth}
		\centering
		\includegraphics[width=\textwidth]{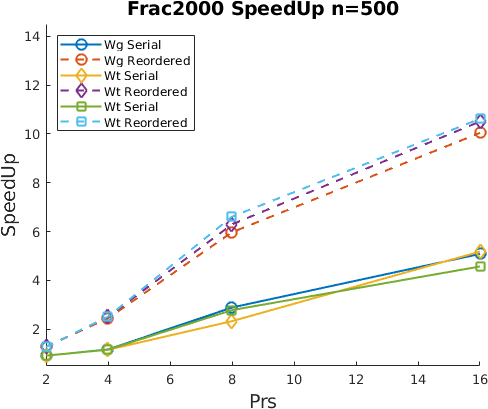}
		\caption{$Frac2000\, n=500$: SpeedUp. }	
	\end{subfigure}
	\begin{subfigure}[t]{.45\textwidth}
		\centering
		\includegraphics[width=\textwidth]{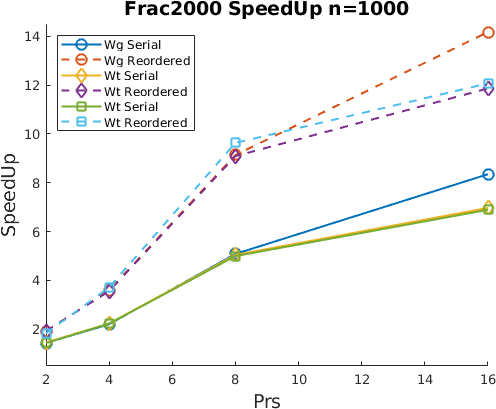}
		\caption{$Frac2000\, n=1000$: SpeedUp. }
		
	\end{subfigure}
	\caption{SpeedUp: $Frac2000$, weighted partitioning strategies base, bipartite and tripartite. Comparison of reordered and serial dofhandler numbering strategies. $n=500$, left, $n=1000$, right.}
	\label{Fig:Frac3000times}
\end{figure*}

In this section we investigate the performances of the solver coupled to
the partitioning strategies introduced in Subsection~\ref{Sub:PartAnalys} (Table~\ref{Notations})
applied to a VEM discretization to the Darcy problem on fractures. In particular we focus on the different performances when the numbering of the Dofs is performed by the simple Algorithm~\ref{Alg:serialDof} and by the Algorithm~\ref{Alg:parDof}. The mesh on which the partition is applied is the minimal mesh \cite{boriodauria} on which we refine uniformly until the number of degrees of freedom satisfies $[\mathcal{F}]\geq\#\mathcal{F}*n$, where $n= 500, 1000$. This approach provides a number of degrees of freedom high enough to justify a parallel approach with few processes. In Figure~\ref{Fig:Frac395meshes} is depicted the final mesh for $Frac512$ with $n=100$.
For this analysis we use the weighted versions of the partitioning strategies aiming at balancing the computational load and minimizing communications among processes. The Dofs are numbered by the basic sequential Algorithm~\ref{Alg:serialDof} and by the reordered Algorithm~\ref{Alg:parDof}, then the resolution times are compared.

Given $t_p$ the resolution time required by $p$ processes, we define the SpeedUp for our parallel resolution as $S_p=\frac{t_1}{t_p}$.

In Figures~\ref{Fig:Frac512times} and~\ref{Fig:Frac3000times} we report the SpeedUp for the cases $Frac512$ and $Frac2000$. As expected the reordered dofhandler displays a clear improved behavior. Figures~\ref{Fig:Frac512times} and~\ref{Fig:Frac3000times} clearly highlights a loss of parallel performances for the case with $16$ processes due to the overloading of shared resources on the CPU. In order to investigate this phenomenon we will present results obtained using different sockets in Subsection~\ref{multisocket}. The tests are performed  on Intel Xeon CPUs with 16 cores \footnote{\url{https://ark.intel.com/content/www/it/it/ark/products/120492/intel-xeon-gold-6130-processor-22m-cache-2-10-ghz.html}}, the $16$ processes tests are performed using all the cores of the CPU. 

Once the efficiency of reordering strategy is proved in terms of SpeedUp, in the following sections the reordered DofHandler is used. The matrix sparsity pattern is presented in Figure~\ref{Figure:Matrices}.

\subsection{Partitioning Analysis}\label{Sub:PartAnalys}

\begin{table}[h]
		\caption{Serial resolution times. 
	}
	\label{Serial}
	\begin{center}
		\footnotesize{
			\begin{tabular}{l|r}%
				\bfseries DFN    & Res. Time(s)  
				\begin{filecontents*}{csvtables/Serial.csv}
Problem ,Processors, NDofs,  Reason Convergence , Residual Norm , Iterations  , Time(sec) ,Partition Cut, Time(sec) , Imbalance, Iterations  , Time(sec) ,Partition Cut, Time(sec) , Imbalance,
512, 	1, 	 	63648 , 2 , 4.88669e-14 ,  2089 , 5.24397 ,  0 , 0  , 1 ,  2089 , 5.24397 ,  0 , 0  , 1 ,
1000,	1,		121424 , 2 , 6.64325e-14 , 3081 , 15.2265 ,  0 , 0  , 1 ,  3081 , 15.2265 ,  0 , 0  , 1 ,
2000,	1, 		151744 , 2 , 8.39113e-14 , 2065 , 13.1338 ,  0 , 0  , 1 ,  2065 , 13.1338 ,  0 , 0  , 1 ,
4000,	1,		332174 , 2 , 6.99009e-14 , 4076 , 57.7411 ,  0 , 0  , 1 ,  4076 , 57.7411 ,  0 , 0  , 1 ,  
\end{filecontents*}
\csvreader[head to column names]{csvtables/Serial.csv}{}
				{\\\hline\csvcoli&\csvcolxii}
			\end{tabular}
		}
	\end{center}
\end{table}
\begin{table}
		\caption{Partition results for graph partitioning $Pg$ and the weighted one $Wg$.  
	}
	\label{TabPgWg_n100}
	\centering
	\footnotesize{
		\begin{tabular}{l|l|r|r|l|l|r|r|r|r}%
			\bfseries DFN  & \bfseries Prs &\multicolumn{2}{c|}{$C$}& \multicolumn{2}{c|}{$I$}  & \multicolumn{2}{c|}{Part. Time(s)} &\multicolumn{2}{c}{Res. Time(s)}  
			\\\hline
			& & $Pg$&$Wg$ & $Pg$&$Wg$& $Pg$&$Wg$& $Pg$&$Wg$ 
			\begin{filecontents*}{csvtables/PartitioningGraph2.csv}
Problem ,Processors, NDofs,  Iterations  , Time(sec) ,Partition Cut, Time(sec) , Imbalance, WEIterations  , WETime(sec) ,WEPartition Cut, WETime(sec) , WEImbalance,
512, 	2, 	 	63648 ,  2091 , 2.676 ,   38 , 0.020 , 0.827 , 2091 , 2.645 ,  50 , 0.022 , 0.970 , 
512, 	4, 	 	63648 ,  2088 , 1.818 ,   90 , 0.023 , 0.642 ,  2085 , 1.403 ,  125 , 0.026 , 0.973 ,  
512, 	8, 	 	63648 ,  2084 , 0.903 , 149 , 0.025 , 0.549 , 2089 , 0.730 , 166 , 0.025 , 0.918 , 
512, 	16,	 	63648 ,  2085 , 0.479 , 360 , 0.054 , 0.534 ,  2082 , 0.420 , 275 , 0.028 , 0.884 ,  
1000,	2, 	 	121424 , 3084 , 8.849 ,  34 , 0.036 , 0.639 ,  3079 , 7.976 ,   45 , 0.038 , 0.944 ,  
1000,	4, 	 	121424 , 3080 , 4.840 , 78 , 0.041 , 0.599 , 	3070 , 3.991 ,  112 , 0.050 , 0.933 , 
1000,	8, 	 	121424 , 3068 , 2.820 , 130 , 0.058 , 0.512 , 	3071 , 2.129 ,   169 , 0.071 , 0.910 , 
1000,	16,	 	121424 , 3066 , 1.372 , 213 , 0.081 , 0.385 ,  3064 , 1.191 , 296 , 0.102 , 0.812 ,
2000,	2, 		151744 , 2060 , 6.726 , 88 , 0.057 , 0.969 ,    2062 , 6.794 , 93 , 0.058 , 0.955 ,
2000,	4, 		151744 , 2057 , 3.482 , 152 , 0.060 , 0.884 ,   2057 , 3.618 , 227 , 0.067 , 0.945 ,  
2000,	8, 		151744 , 2055 , 1.979 , 293 , 0.067 , 0.760 ,   2057 , 1.857 , 339 , 0.068 , 0.941 ,  
2000,	16,		151744 , 2054 , 1.092 , 455 , 0.084 , 0.608 ,  2055 , 1.000 , 551 , 0.076 , 0.890 , 
4000,	2, 		332174 , 3785 , 28.971 , 86 , 0.118 , 0.901 ,   3776 , 28.583 , 101 , 0.155 , 0.965 ,  
4000,	4, 		332174 , 3762 , 15.022 , 177 , 0.120 , 0.827 ,  3758 , 14.141 , 204 , 0.128 , 0.970 , 
4000,	8, 		332174 , 3757 , 8.200 , 372 , 0.136 , 0.731 ,  3758 , 7.608 , 501 , 0.144 , 0.930 , 
4000,	16,		332174 , 3751 , 5.427 , 556 , 0.137 , 0.537 ,  3756 , 4.629 , 652 , 0.141 , 0.950 ,  
\end{filecontents*}
\csvreader[head to column names]{csvtables/PartitioningGraph2.csv}{}
			{\\\hline\csvcoli&\csvcolii&\csvcolvi&\csvcolxi&\csvcolviii&\csvcolxiii&\csvcolvii&\csvcolxii&\csvcolv&\csvcolx}
		\end{tabular}
	}
\end{table}
\begin{table}
	\caption{Partition results for bipartite partitioning $Pb$ and the weighted one $Wb$.
	}
	\label{TabPbWb_n100}
	\centering
	\footnotesize{
		
		\begin{tabular}{l|l|r|r|l|l|r|r|r|r}%
			\bfseries DFN  & \bfseries Prs &\multicolumn{2}{c|}{$C$}& \multicolumn{2}{c|}{$I$}  & \multicolumn{2}{c|}{Part. Time(s)} & \multicolumn{2}{c}{Res. Time(s)}  
			\\\hline
			& & $Pb$&$Wb$ & $Pb$&$Wb$& $Pb$&$Wb$ & $Pb$&$Wb$
			\begin{filecontents*}{csvtables/PartitioningBipartite2.csv}
Problem ,Processors, NDofs,  Iterations  , Time(sec) ,Partition Cut, Time(sec) , Imbalance, WEIterations  , WETime(sec) ,WEPartition Cut, WETime(sec) , WEImbalance,
512, 	2, 	 	63648 , 2090 , 2.682 , 40 , 0.020 , 0.864 ,    2088 , 2.648 ,  60 , 0.022 , 0.988 ,  
512, 	4, 	 	63648 , 2087 , 1.612 , 80 , 0.022 , 0.783 ,     2086 , 1.412 ,  151 , 0.028 , 0.937 , 
512, 	8, 	 	63648 , 2089 , 0.816 , 145 , 0.024 , 0.701 ,  2082 , 0.748 , 219 , 0.030 , 0.887 ,   
512, 	16,	 	63648 , 2085 , 0.474 , 211 , 0.026 , 0.500 ,  2083 , 0.452 , 290 , 0.032 , 0.836 ,
1000,	2, 	 	121424 , 3078 , 8.122 , 45 , 0.039 , 0.888 ,   3077 , 7.959 , 50 , 0.039 , 0.980 ,  
1000,	4, 	 	121424 ,  3072 , 4.090 ,90 , 0.041 , 0.786 ,  3072 , 4.001 , 102 , 0.053 , 0.943 ,  
1000,	8, 	 	121424 , 3068 , 2.389 , 140 , 0.053 , 0.711 ,  3072 , 2.144 , 180 , 0.052 , 0.896 ,  
1000,	16,	 	121424 , 3069 , 1.399 , 239 , 0.062 , 0.553 ,   3070 , 1.245 , 315 , 0.070 , 0.863 , 
2000,	2, 		151744 , 2059 , 7.006 , 76 , 0.056 , 0.967 ,   2061 , 6.782 , 93 , 0.057 , 0.962 ,  
2000,	4, 		151744 ,  2056 , 3.336 , 157 , 0.061 , 0.954 ,  2058 , 3.379 , 226 , 0.065 , 0.952 , 
2000,	8, 		151744 , 2059 , 1.821 , 335 , 0.063 , 0.775 ,   2058 , 1.837 , 420 , 0.073 , 0.914 ,  
2000,	16,		151744 , 2055 , 1.060 , 497 , 0.079 , 0.743 , 2055 , 1.008 , 611 , 0.080 , 0.867 ,  
4000,	2, 		332174 , 3774 , 28.428 , 83 , 0.118 , 0.992 ,  3793 , 28.045 , 92 , 0.119 , 0.986 ,  
4000,	4, 		332174 , 3762 , 14.781 , 176 , 0.121 , 0.859 ,  3759 , 14.699 , 195 , 0.124 , 0.936 ,
4000,	8, 		332174 , 3751 , 7.512  ,  348 , 0.130 , 0.820 , 3758 , 7.303 , 494 , 0.138 , 0.929 ,  
4000,	16,		332174 , 3754 , 5.070 , 603 , 0.135 , 0.676 ,   3754 , 4.663 ,  767 , 0.155 , 0.902 , 
\end{filecontents*}
\csvreader[head to column names]{csvtables/PartitioningBipartite2.csv}{}
			{\\\hline\csvcoli&\csvcolii&\csvcolvi&\csvcolxi&\csvcolviii&\csvcolxiii&\csvcolvii&\csvcolxii&\csvcolv&\csvcolx}
		\end{tabular}
	}
\end{table}
\begin{table}
	\caption{Partition results for tripartite partitioning $Pt$ and the weighted one $Wt$.}
\label{TabPtWt_n100}
	\centering
	\footnotesize{
		\begin{tabular}{l|l|r|r|l|l|r|r|r|r}%
			\bfseries DFN  & \bfseries Prs &\multicolumn{2}{c|}{$C$}& \multicolumn{2}{c|}{$I$}  & \multicolumn{2}{c|}{Part. Time(s)}  &\multicolumn{2}{c}{Res. Time(s)}  
			\\\hline
			& & $Pt$&$Wt$ & $Pt$&$Wt$& $Pt$&$Wt$& $Pt$&$Wt$ 
			\begin{filecontents*}{csvtables/PartitioningTripartite2.csv}
Problem ,Processors, NDofs,  Iterations  , Time(sec) ,Partition Cut, Time(sec) , Imbalance, WEIterations  , WETime(sec) ,WEPartition Cut, WETime(sec) , WEImbalance,
512, 	2, 	 	63648 , 3082 , 2.752 ,   50 , 0.021 ,  0.872 ,  3075 , 2.645,  70 , 0.024 , 0.986 ,  
512, 	4, 	 	63648 , 2090 , 1.579 ,  93 , 0.022 ,  0.801 ,  2088 , 1.368 ,  143 , 0.026 , 0.963 ,  
512, 	8, 	 	63648 , 2089 , 0.816 , 136 , 0.025 , 0.864 ,  2088 , 0.741 , 136 , 0.025 , 0.864 ,   
512, 	16,	 	63648 , 2083 , 0.478 , 219 , 0.027 , 0.466 ,  2083 , 0.452 , 290 , 0.032 , 0.836 , 
1000,	2, 	 	121424 , 3075 , 8.004 , 45 , 0.037 , 0.886 ,    3082 , 8.192 ,  45 , 0.039 , 0.961 , 
1000,	4, 	 	121424 , 3077 , 4.448 , 87 , 0.040 , 0.805 ,	3076 , 4.018 , 112 , 0.057 , 0.947 ,  
1000,	8, 	 	121424 , 3077 , 2.402 , 140 , 0.053 , 0.760 ,   3071 , 2.203 , 189 , 0.064 , 0.896 ,  
1000,	16,	 	121424 , 3068 , 1.360 , 232 , 0.063 , 0.594 ,  3065 , 1.132 , 308 , 0.058 , 0.863 ,
2000,	2, 		151744 , 2060 , 6.628 , 79 , 0.055 , 0.961 ,  2062 , 6.922 ,  165 , 0.065 , 0.974 , 
2000,	4, 		151744 , 2059 , 3.361 ,  166 , 0.062 , 0.924 ,  2059 , 3.501 , 192 , 0.064 , 0.958 ,  
2000,	8, 		151744 , 2059 , 2.029 , 312 , 0.066 , 0.796 ,   2057 , 1.827 ,  356 , 0.069 , 0.929 ,  
2000,	16,		151744 , 2055 , 1.091 , 472 , 0.073 , 0.761 , 2055 , 0.962 , 547 , 0.080 , 0.919 ,
4000,	2, 		332174 , 4073 , 30.125 , 83 , 0.117 , 0.973 ,  4078 , 31.120 , 98 , 0.118 , 0.979 , 
4000,	4, 		332174 , 3763 , 14.415 , 175 , 0.121 , 0.934 , 3762 , 13.857 , 282 , 0.134 , 0.951 ,
4000,	8, 		332174 , 3758 , 7.543 , 358 , 0.129 , 0.826 ,  3758 , 7.215 , 453 , 0.138 , 0.927 , 
4000,	16,		332174 , 3753 , 5.225 , 577 , 0.138 , 0.743 , 3751 , 5.106 , 728 , 0.149 , 0.880 ,  
\end{filecontents*}
\csvreader[head to column names]{csvtables/PartitioningTripartite2.csv}{}
			{\\\hline\csvcoli&\csvcolii&\csvcolvi&\csvcolxi&\csvcolviii&\csvcolxiii&\csvcolvii&\csvcolxii&\csvcolv&\csvcolx}
		\end{tabular}
	}
\end{table}
In this section we compare the performances of the six partitioning strategies presented in Section~\ref{Sec:ParallelPartitioning} for different DFNs on a mesh with approximately $[\mathcal{F}]\geq \#\mathcal{F} *n,$ with $n=100$ degrees of freedom. Here we only consider the reordered version of the Dofs indices.

In Table~\ref{Serial} we report the resolution times for the serial problems. In Tables~\ref{TabPgWg_n100}, \ref{TabPbWb_n100} and \ref{TabPtWt_n100} we report the results comparing the three partitioning strategies in their unweighted and weighted versions. In the tables the second column reports the number of processes involved for each test case specified on first column. The columns $3$ and $4$ report the cut $C$ for unweighted and weighted partitioning, respectively. The columns $5$ and $6$ the imbalance $I$.  $Part.Time(s)$ on columns $7$ and $8$ is the partitioning time employed by METIS, the last two columns report the resolution time employed by the PETSc.

\begin{figure*}[h]
	\begin{subfigure}[t]{.45\textwidth}
		\centering
		\includegraphics[width=\textwidth]{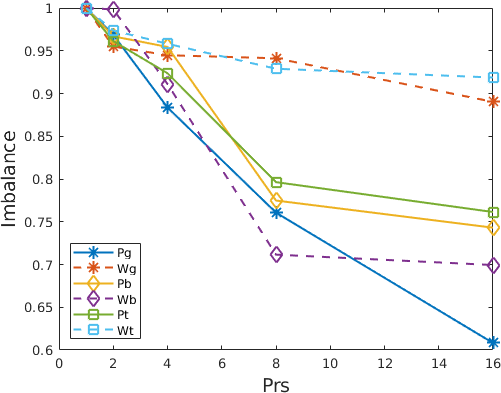}
	\end{subfigure}
	~
	\begin{subfigure}[t]{.45\textwidth}
		\centering
		\includegraphics[width=\textwidth]{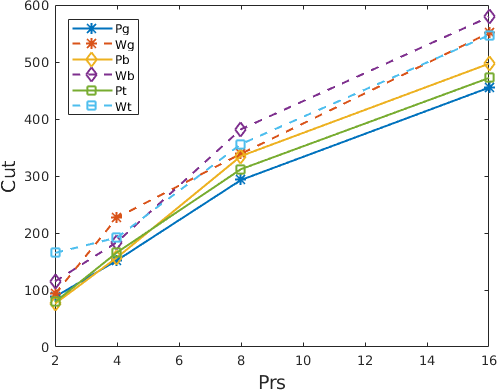}
	\end{subfigure}
	\caption{$Frac2000$: imbalance (left) and cut (right) increasing the number of processes, $n=100$.}
	\label{Figure:2000Imbalance}
\end{figure*}
\begin{figure*}
	\includegraphics[width=.4\textwidth]{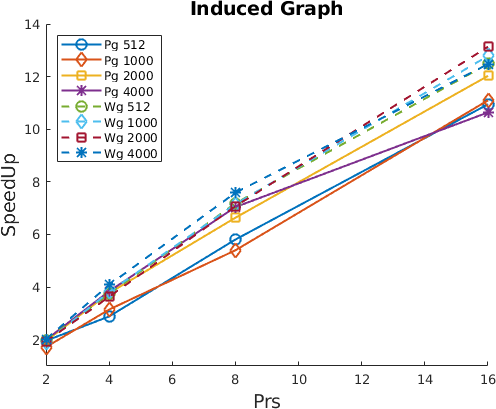}%
	~
	\includegraphics[width=.4\textwidth]{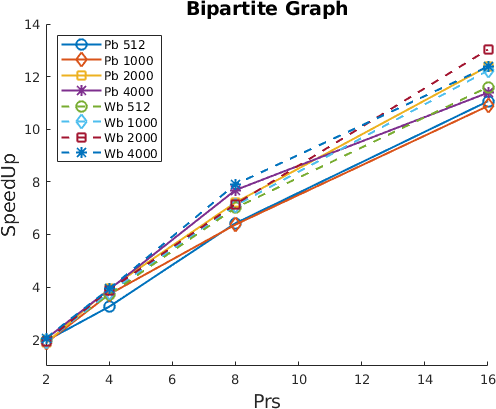}%
	~
	
	\includegraphics[width=.4\textwidth]{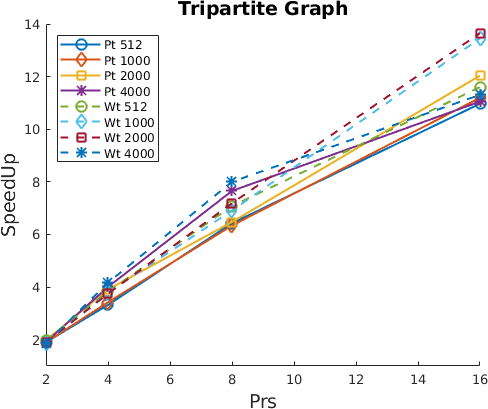}%
	\caption{SpeedUp comparison among partitioning strategies. The average number of degrees of freedom on each fracture is $n=100$.}
	\label{Figure:BipartiteSpeedUp}
\end{figure*}

\begin{figure*}
	\includegraphics[width=.45\textwidth]{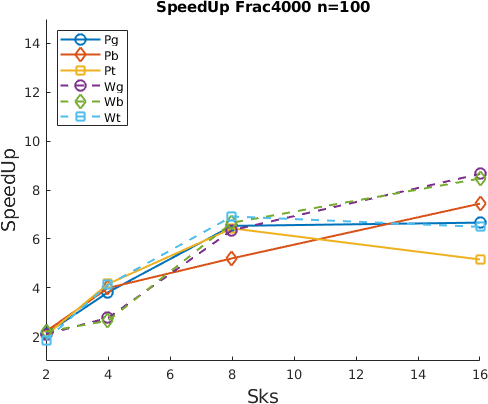}%
	\includegraphics[width=.45\textwidth]{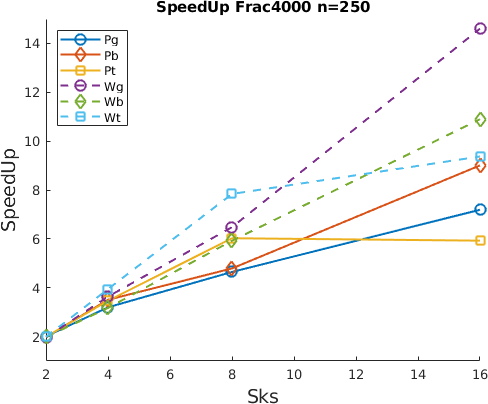}%
	~
	
	\includegraphics[width=.45\textwidth]{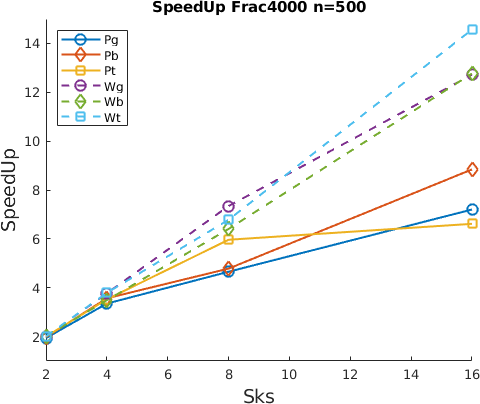}	
	\includegraphics[width=.45\textwidth]{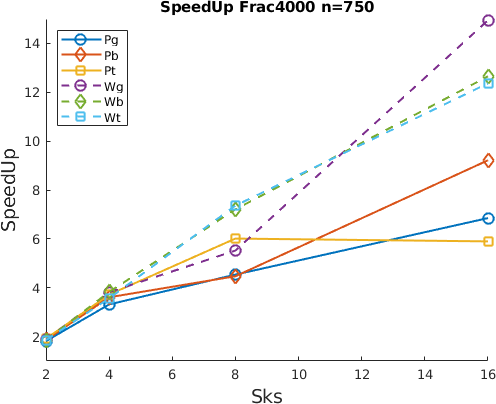}
	\caption{Multi-socket SpeedUp for $Frac4000$,   $[\mathcal{F}]$: $332174$, $1260833$ , $1889468$ and $2789051$. The average number of degrees of freedom on each fracture is $n$; partitioning strategies comparison for their weighted and unweighted versions.}
	\label{Figure:MultiNodesT}
\end{figure*}

For each partitioning strategy, comparing the results of the unweighted and weighted versions we can see that the cut $C$ is increasing for the weighted versions because the partitioner aims at minimizing the number of degrees of freedom on the cut traces instead of the number of the cut traces. The cut $C$ represents the number of communicating traces of the DFN, that are edges of the graph representing the DFN only for the $Pg$ and $Wg$ strategies. In the weighted version the cuts concentrate on more shorter traces. We note that an increased number of communicating traces arises in the bipartite and tripartite case, as long as they focus on the data load and not on the quantity of the cut edges.

Concerning the imbalance $I$, we see that the weighted versions have an imbalance closer to one and less dependent on the number of processes, that denotes a more uniform distribution of the workload among the processes.
In Figure~\ref{Figure:2000Imbalance} we report for $Frac2000$ plots of imbalance for all the partitioning strategies, the advantages in terms of Dofs imbalance is evident for base graph and tripartite graph in front of a negligible loss for the cut $C$.

From columns $7-8$ ($Part. Time(s)$) we can see that the partitioning time is almost equal for all the partitioning strategies.

In Figure~\ref{Figure:BipartiteSpeedUp} we compare the SpeedUp for the weighted and non-weighted strategies, we can observe that they are quite similar and the weighted version behaves slightly better. In this particular case the total number of Dofs is not so large to overload the shared resources of the CPU and the degradation seen in Figure~\ref{Fig:Frac3000times} is slightly appreciable for the $Frac4000$.

\subsection{Multi-Sockets Analysis}
\label{multisocket}

As already noted in the previous sections when the number of partitions is higher than $8$ and the number of Dofs is quite large our tests clearly highlight a loss of SpeedUp. This phenomenon can be attributed to an overloading of the common resources of the CPUs. In order to confirm this interpretation we repeat the same tests entrusting only a process to each CPU.
Being the communications between different sockets more expensive \cite{balay2019petsc}  we focus on a larger problem that is the $Frac4000$  on which we increase the number of degrees of freedom $[\mathcal{F}]\geq\#\mathcal{F}*n$ from $n= 100$ to $n=250$, $n=500$ and $n=750$. In Figure~\ref{Figure:MultiNodesT} we report the SpeedUp results, that highlight better performances increasing the number of degrees of freedom due to the increased workload of the processes with respect to the cost of communications. In this case the degeneration of performances passing from 8 to 16 processes is less relevant because there is not competition between processes in the use of the shared resources of the CPU.

\subsection{Mesh Partitioning}\label{Sub:MP}

In this section we consider the standard partitioning of the Dofs based on the connectivity of the mesh. In particular we construct the adjacency matrix of the graph of the degrees of freedom, i.e., we have one graph-node for each mesh-node (Dof) and one graph-edge for each mesh-edge. We apply a partitioning of this graph using a cut-edge minimization strategy; with this approach the amount of data communicated at each iteration among the processes should be minimized. 

\begin{figure*}[!h]
	\begin{subfigure}[t]{.4\textwidth}	\includegraphics[width=\textwidth]{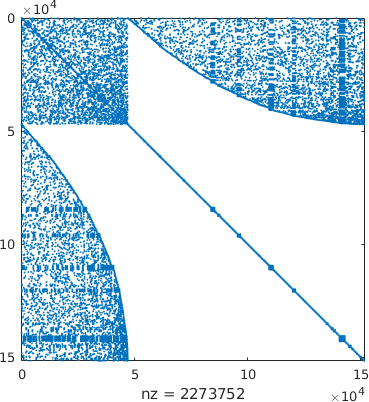}
		\caption{Serial DofHandler matrix structure. }
		\label{Figure:Adj}
	\end{subfigure}
	\begin{subfigure}[t]{.4\textwidth}	\includegraphics[width=\textwidth]{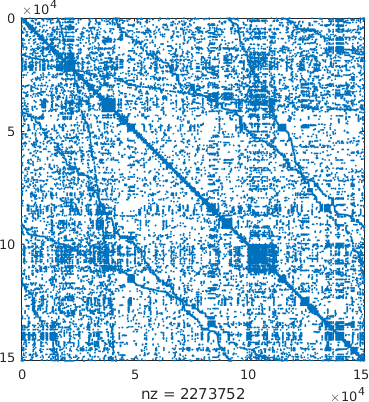}
	\caption{$MeshP$ partitioning. }
	\label{Figure:AdjMesh}
\end{subfigure}
	\begin{subfigure}[t]{.4\textwidth}	\includegraphics[width=\textwidth]{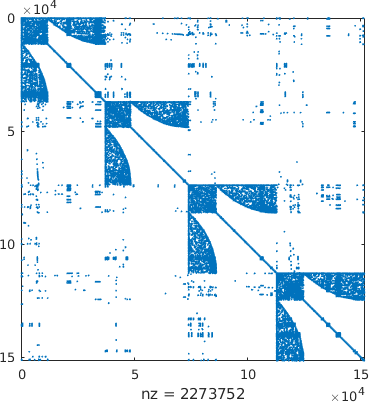}
	\caption{$Wg$ partitioning. }
	\label{Figure:AdjBase}
\end{subfigure}
	\begin{subfigure}[t]{.4\textwidth}	\includegraphics[width=\textwidth]{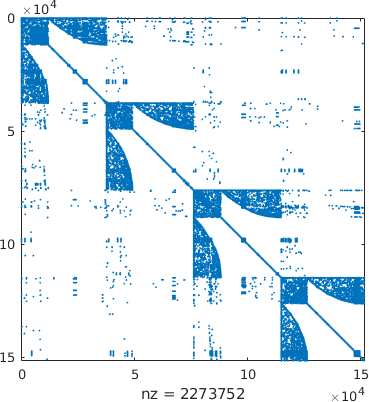}
	\caption{$Wb$ partitioning. }
	\label{Figure:AdjBip}
\end{subfigure}
	\begin{subfigure}[t]{.4\textwidth}	\includegraphics[width=\textwidth]{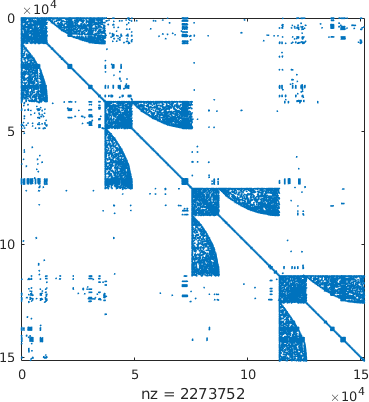}
	\caption{$Wt$ partitioning. }
	\label{Figure:AdjTrip}
\end{subfigure}
	\caption{$Frac2000$: Sparsity pattern. Partitioning among 4 processes}
	\label{Figure:Matrices}
\end{figure*}

The partitioning of the Dofs is handled by METIS. Then the resulting set of indices on each process is numbered in a contiguous way (see Figures~\ref{Figure:Adj}, \ref{Figure:AdjMesh}).

\begin{figure*}[h]
	\begin{subfigure}[t]{.45\textwidth}
		\includegraphics[width=\textwidth]{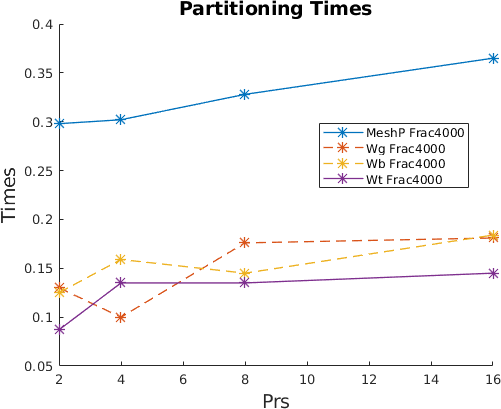}
		\caption{$Frac4000$: Partitioning times employed by METIS comparison.}
		\label{Figure:PartTimes}
	\end{subfigure}
	\begin{subfigure}[t]{.45\textwidth}	\includegraphics[width=\textwidth]{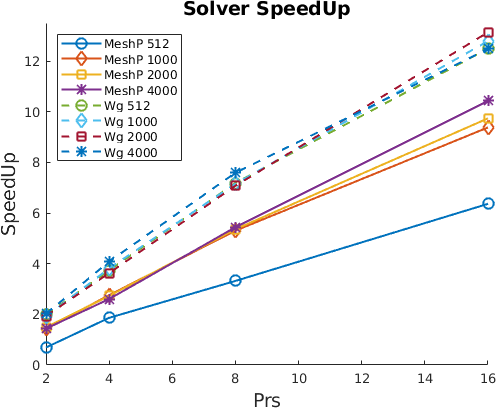}
		\caption{Mesh partition SpeedUp. Comparison with weighted graph partition. }
		\label{Figure:ScalMeshPart}
	\end{subfigure}
	\caption{Induced mesh partitioning results for partitioning times and solver SpeedUp.}
	\label{Figure:MeshPartRes}
\end{figure*}

The partitioning time employed by METIS for this approach is much higher with respect to the ones of the previously presented partitioning strategies due to the larger graph to be partitioned, see Figure~\ref{Figure:PartTimes}. Nevertheless, the solver SpeedUp is quite similar, see Figure~\ref{Figure:ScalMeshPart}.  When the number of processes increases the DFN-based graph partitionings present a clear better overall behavior, \cite{ushijima2019multilevel}. 
In Table~\ref{TabMesh} we report the partitioning time and the resolution time for the same mesh ($n=100$) used for the test cases of Tables~\ref {TabPgWg_n100}, \ref{TabPbWb_n100} and \ref{TabPtWt_n100}. We can observe that the resolution time for the mesh partitioning is always larger with respect to the proposed partitioning strategies.
This behavior can be explained observing the sparsity pattern of the matrix. In Figure~\ref{Figure:Matrices} we report the spy for the serial DofHandler, for the mesh partitioning and for the weighted graph partitionings proposed. The large number of off-diagonal block elements of the mesh partitioning approach is responsible for a lower efficiency of the matrix vector products, whereas in the reordered weighted versions we have similar structures with marginal differences and a higher clustering of the nonzeros elements on the rows that helps improving the matrix vector product.

\begin{table}[h]
		\caption{Mesh Partitioning.}
	\label{TabMesh}
	\begin{center}
		\footnotesize{
			\begin{tabular}{l|l|r|r}%
				\bfseries DFN  & \bfseries Prs  & Part. Time(s) & Res.Time(s) 
				
				\begin{filecontents*}{csvtables/MetisDatiTable2.csv}
 1,	63648 , 512 , 4.78261e-14 , 2090 , 0 , 	5.025 , 44.0799 , 44.0799,	0 , 	 , 0 	 , 63648 	 , 1	 , 0 	 , 
1,	63648 , 512 , 4.78261e-14 , 2090 , 0 , 	5.025 , 44.0799 , 44.0799,	0 , 	 , 0 	 , 63648 	 , 1	 , 0 	 , 
 2,	63648 , 512 , 4.92653e-14 , 2089 , 0 , 	3.564 , 52.6908 , 52.6908, 0.052 , 	 , 0 	 , 31824 	 , 1	 , 107432 	 , 	
 4,	63648 , 512 , 4.83389e-14 , 2086 , 0 , 	1.836 , 52.8862 , 52.8862, 0.055	 , 	 , 0 	 , 15912 	 , 1	 , 166986 	 , 	
 8,	63648 , 512 , 4.84169e-14 , 2085 , 0 , 	1.010 , 52.7041 , 52.7041, 0.060 , 	 , 0 	 , 7956 	 , 1	 , 199256 	 , 
 16,	63648 , 512 , 4.87095e-14 , 2085 , 0 , 	0.570 , 53.2615 , 53.2615, 0.068 , 	 , 0 	 , 3978 	 , 1	 , 225617 	 , 
 1,	121424 ,1000 , 6.31545e-14 ,3085 ,0 ,	14.789 ,82.2862 ,82.2862, 0	 , 	 , 0 	 , 121424 	 , 1	 , 0 	 , 	
 2,	121424 ,1000 , 6.63071e-14 ,3080 ,0 ,	10.239 ,98.5856 ,98.5856, 0.102	 , 	 , 0 	 , 60712 	 , 1	 , 203722 	 , 	
 4,	121424 ,1000 , 6.57599e-14 ,3069 ,0 ,	5.363 ,98.8478 ,98.8478, 0.105	 , 	 , 0 	 , 30356 	 , 1	 , 323698 	 , 	
 8,	121424 ,1000 , 6.60105e-14 ,3071 ,0 ,	2.789 ,99.1687 ,99.1687, 0.112	 , 	 , 0 	 , 15178 	 , 1	 , 376901 	 , 	
 16,	121424 ,1000 , 6.29453e-14 ,3065 ,0 ,	1.577 ,99.6427 ,99.6427, 0.135	 , 	 , 0 	 , 7589 	 , 1	 , 429736 	 , 	
1,	151744 ,2000 , 8.10877e-14 ,2063 ,0 ,	12.648 ,99.2588 ,99.2588, 0	 , 	 , 0 	 , 151744 	 , 1	 , 0 	 , 	
2,	151744 ,2000 , 8.51657e-14 ,2061 ,0 ,	8.517 ,123.384 ,123.384, 0.133	 , 	 , 0 	 , 75872 	 , 1	 , 267414 	 , 	
4,	151744 ,2000 , 8.33008e-14 ,2059 ,0 ,	4.638 ,127.935 ,127.935, 0.146	 , 	 , 0 	 , 37936 	 , 1	 , 466003 	 , 
8,	151744 ,2000 , 8.46873e-14 ,2056 ,0 ,	2.360 ,119.388 ,119.388, 0.146	 , 	 , 0 	 , 18968 	 , 1	 , 506598 	 , 
16,	151744 ,2000 , 8.55163e-14 ,2053 ,0 ,	1.299 ,124.12 ,124.12,	0.170	 , 	 , 0 	 , 9484 	 , 1	 , 562416 	 , 
1,	332174 ,4000 , 7.33097e-14 ,4080 ,0 ,	56.103 ,216.51 ,216.51,	0	 , 	 , 0 	 , 332174 	 , 1	 , 0 	 , 
2,	332174 ,4000 , 7.41575e-14 ,4076 ,0 ,	39.084 ,261.617 ,261.617, 0.298	 , 	 , 0 	 , 166087 	 , 1	 , 575876 	 , 
4,	332174 ,4000 , 7.35789e-14 ,4065 ,0 ,	21.548 ,262.883 ,262.883, 0.302	 , 	 , 0 	 , 83044 	 , 0.999988	 , 1011162 	 , 	
8,	332174 ,4000 , 7.39274e-14 ,3763 ,0 ,	10.316 ,263.914 ,263.914, 0.328	 , 	 , 0 	 , 41522 	 , 0.999976	 , 1094796 	 , 	
16,	332174 ,4000 , 7.40355e-14 ,3755 ,0 ,	5.382 ,262.476 ,262.476, 0.365	 , 	 , 0 	 , 20761 	 , 0.999952	 , 1199541 	 , 	
\end{filecontents*}
\csvreader[head to column names]{csvtables/MetisDatiTable2.csv}{}
				{\\\hline
					\csvcoliii &\csvcoli &  \csvcolx & \csvcolvii}
				
			\end{tabular}
		}
	\end{center}
\end{table}

\subsection{An adaptive VEM mesh refinement test}\label{Sub:VEM}

\begin{figure*}
	\includegraphics[width=.5\textwidth]{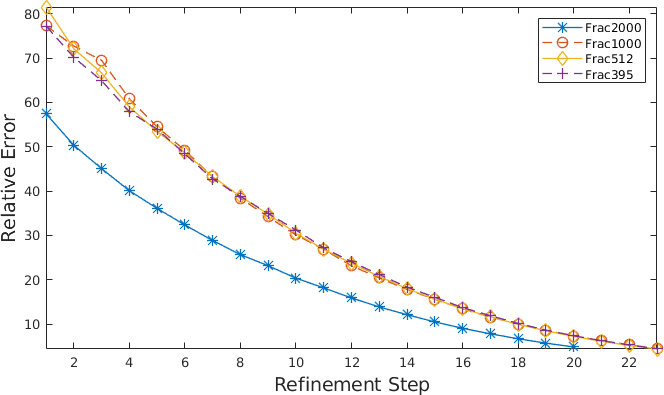}%
	\includegraphics[width=.5\textwidth]{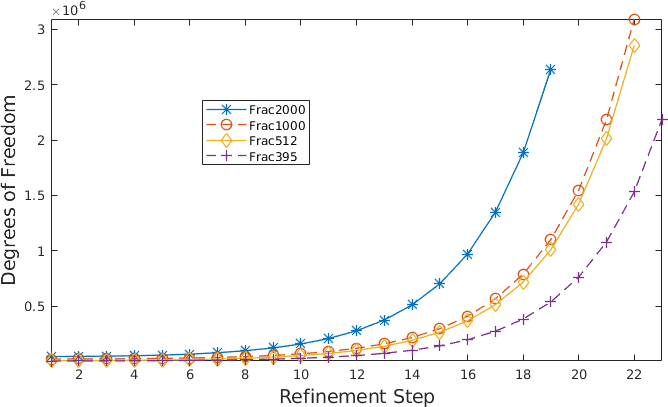}
	\caption{Relative error (left) and increasing degrees of freedom (right) during an adaptive mesh refinement.}
	\label{Fig:VEMCurves}
\end{figure*}
\begin{figure}
	\begin{center}
		\includegraphics[width=.4\textwidth]{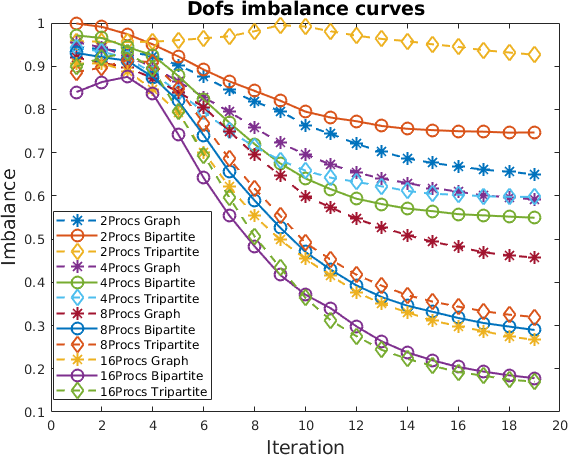}
	\end{center}
	\caption{$Frac2000$: imbalance curves during VEM resolution. }
	\label{Fig:ImbalVEM2000}
\end{figure}

In this last section we investigate the interplay between partitioning and mesh refinement, and, in particular, the degeneration of the performances of the proposed partitioning strategies when the partitioning is performed on a starting mesh that is subject to an iterative mesh refinement based on {\sl a posteriori} error estimates, \cite{Berrone-Borio:2017}. In Figure~\ref{Fig:VEMCurves} the convergence curves of the relative error are reported on the left; on the right the growth of Dofs number due to the mesh refinement, \cite{boriodauria}.
As long as it is not possible to provide the adaptive refinement \emph{a priori}, a computationally costly, but performing approach would apply a partitioning of the DFN during each refinement step. Another approach is to partition the DFN every fixed predefined refinement steps (for example each five refinements). However we aim at observing the behavior of our partitioning strategies without a re-partitioning phase to see how rapidly the initial partitioning degenerate.
We apply the partitioning strategies on the minimal mesh at the beginning of the computation. We do not expect, as long as the refining is not uniform, that the imbalance can remain constant during the resolution. Curves presented in Figure~\ref{Fig:ImbalVEM2000} suggest to apply a re-partitioning each five refinement iterations.

\section{Conclusions}
\label{Sec:Concl}
In this article we present three partitioning strategies for the DFN flow simulations. For each of them, we analyze several performances parameters, among them the cut, the imbalance, the partitioning time and the resolution time. The weighted partitioning strategies in general perform better, moreover they have similar behavior, with in general slight loss in performances for the bipartite graph partitioning. The proposed partitioning strategies are computationally  cheaper than the classical mesh induced one in partitioning time and perform better in resolution time. 

The reordering aims at avoiding that one or few processes are overloaded of communications with respect to the other processes, trying to equally distribute the communication among all the processes as presented in Section~\ref{loadDofs}.

For problems with a large number of degrees of freedom we also compare the SpeedUp on processes running on the same socket and on different CPUs clearly highlighting a degradation of the SpeedUp when the CPU and memory access resources are overloaded. The proposed methods preserve the SpeedUp with a higher dimension of the DFN problem.

\subsubsection*{Acknowledgements}
This work is supported by the MIUR project ``Dipartimenti di Eccellenza 2018-2022'' (CUP E11G18000350001), PRIN project ``Virtual Element Methods: Analysis and Applications'' (201744KLJL\_004) and by INdAM-GNCS. Computational resources supported by HPC@polito.it and SmartData@polito.


\bibliographystyle{plain}       
{\footnotesize
\bibliography{bibliography}}   

\end{document}